\documentclass[review,3p]{elsarticle}

\usepackage[T1]{fontenc}

\usepackage[english]{babel}
\usepackage{amsmath}
\usepackage{dsfont}
\usepackage{graphicx}
\usepackage{hyperref}
\hypersetup{colorlinks=true, allcolors=blue}
\usepackage{enumitem}
\usepackage{makecell}
\usepackage{bbold}
\usepackage{commath}

\usepackage{standalone}
\usepackage{tikz} 
\usepackage{xfp}
\usetikzlibrary{shapes,arrows}

\usepackage[labelsep=period]{caption}
\usepackage{multirow}
\usepackage{multicol}

\usepackage{lineno}

\usepackage{bm}
\usepackage{comment}
\usepackage{pbox}
\usepackage{xcolor}
\usepackage{subfigure}

\graphicspath{ {image/} }
\journal{Structural and Multidisciplinary Optimization}

\usepackage{lipsum}
\makeatletter \def\ps@pprintTitle{\let\@oddhead\@empty \let\@evenhead\@empty \def\@oddfoot{\hfill\today} \let\@evenfoot\@oddfoot} \makeatother

\begin{document}

\begin{frontmatter}

\title{S-BORM: Reliability-based optimization of general systems using buffered optimization and reliability method}

\author[1]{Ji-Eun Byun\corref{cor1}}
\ead{Ji-Eun.Byun@glasgow.ac.uk}
\address[1]{James Watt School of Engineering, University of Glasgow, Glasgow, United Kingdom}

\author[2]{Welington de Oliveira}
\address[2]{Mines Paris, Universit\'e PSL, Centre de Math\'ematiques Appliqu\'ees (CMA), Sophia Antipolis, France}

\author[3]{Johannes O. Royset}
\address[3]{Operations Research Department, Naval Postgraduate School, California, United States}

\cortext[cor1]{Corresponding author}

\begin{abstract}
    Reliability-based optimization (RBO) is crucial for identifying optimal risk-informed decisions for designing and operating engineering systems. However, its computation remains challenging as it requires a concurrent task of optimization and reliability analysis. Moreover, computation becomes even more complicated when considering performance of a general system, whose failure event is represented as a link-set of cut-sets. This is because even when component events have smooth and convex limit-state functions, the system limit-state function has neither property, except in trivial cases. To address the challenge, this study develops an efficient algorithm to solve RBO problems of general system events. We employ the buffered optimization and reliability method (BORM), which utilizes, instead of the conventional failure probability definition, the buffered failure probability. The proposed algorithm solves a sequence of difference-of-convex RBO models iteratively by employing a proximal bundle method. For demonstration, we design three numerical examples with increasing complexity that includes up to 108 cut-sets, which are solved by the proposed algorithm within a minute with high accuracy. We also demonstrate its robustness by performing extensive parametric studies. 
\end{abstract}


\end{frontmatter}

\section{Introduction}

To secure disaster resilience of a community, it is crucial to make optimal decisions when designing and operating engineering systems (e.g., structural systems, infrastructures, and mechanical systems) and these decisions should account appropriately for hazard risks. This can be done by performing reliability-based optimization (RBO), where a design cost is minimized while satisfying reliability constraints \cite{rbo94, mbnDm20}. A common, reasonable way to define such reliability constraints is to constrain a failure probability under a target level. However, such a combined task of probabilistic analysis and optimization makes RBO problems theoretically and computationally challenging. Moreover, RBO problems become even more challenging when an event of interest is a system event, whose performance is determined by joint performance of multiple component events. Nonetheless, to enable accurate decision-making, it is critical to consider interdependent component events simultaneously \cite{bns_jt21,sra12}.

To formulate an RBO problem, a system failure probability needs to be represented as a function of design variables. Since explicit expression of such functions is often unavailable, various approximate methods have been proposed. This includes replacing failure probability by reliability index \cite{sl-rbo11, rbo-decoup01} or failure rate \cite{inm_pipe22}, evaluating probabilities approximately using surrogate models such as kriging \cite{rbo-kriging21,rbo_surro_mcs19} or co-kriging \cite{Perdikaris}, simplifying optimization 
by approximate functions \cite{mbnDm20,rbo_polytop19}, and leveraging metaheuristic optimization such as genetic algorithms \cite{rbo-ga21} or particle swarm optimization algorithms \cite{rbo-pso22}.

While these previous studies provide characteristic advantages of their own, another promising approach is to directly utilize realizations of random variables, i.e., samples or data points. This strategy exempts us from deriving problem-specific formulas. It enables data-driven optimization, which is particularly favored when underlying parametric distributions of given data are unknown. However, the definition of failure probability makes it challenging to handle realizations during optimization. This is because a sample is assigned either 1 or 0 depending on whether it lies in the failure domain or not. Such binary assignment lacks gradient information. This makes it difficult to employ efficient gradient-based optimization algorithms \cite{optim_primer}. Although \cite{rbo_samp07} proposed formulations to calculate such gradients, some restrictive conditions are necessary (e.g., random variables follow the multivariate normal distribution or transformation thereof). The challenge can be addressed by utilizing, instead of the conventional failure probability, the {\it buffered failure probability} \cite{bpf10, borm22, bpf_grad21, rbdo_cvar22}. This way of defining reliability permits (sub)gradient information even in a data-driven setting \cite{bpf_grad21}. Motivated by this fact, \cite{borm22} proposed the {\it buffered optimization and reliability method} (BORM), which enables efficient data-driven optimization of reliability. It is noted that the use of the buffered failure probability still closely aligns with risk management by the conventional failure probability as the two failure probabilities have strongly positive correlations \cite{borm22}. 

General systems are often represented as a series system of multiple parallel systems or equivalently, a link-set of cut-sets. This representation is useful since it can cover any system types. Howver, it also highlights the challenges associated with optimization of general systems: A system event tends to have a nonsmooth and nonconvex limit-state function with respect to design variables even when component events have smooth and convex limit-state functions \cite{bpf10}. This effect is unavoidable, but is made more manageable in BORM because it avoids the 0-1 discontinuity caused by failure probability formulations.

Motivated by the promising properties of BORM, this study proposes a novel algorithm for reliability-based optimization, particularly intended to handle general systems. Since the algorithm addresses {\em system} reliability optimization using BORM, it is named {\it S-BORM}. The S-BORM algorithm solves subproblems that are formulated by BORM and uses data or sample points to estimate the buffered failure probability. The subproblems fall into the class of {difference-of-convex} optimization problems
\cite{dc_textbook16} and are handled by the difference-of-convex bundle method of \cite{pb17}.
As we show in the following discussions, the difference-of-convex representation of subproblems can be derived by adaptively linearizing limit-state functions at a current solution candidate. 
We further enhance computational efficiency by employing an active-set strategy, which significantly reduces the number of samples (or data points) that need to be considered in each iteration. 
S-BORM has rigorous convergence properties: for limit-state functions that are linear in the decision variables, the algorithm can only converge to critical points.  We design three numerical examples with increasing complexity to demonstrate that the algorithm works practically well for both linear {\em and} nonlinear limit-state functions. In addition, we perform extensive parametric test, which shows the robustness of the algorithm. The S-BORM algorithm is also developed as a Matlab-based function applicable for customized problems, which is available at \url{https://github.com/jieunbyun/sborm}.

The paper is organized as follows. Section \ref{sec:background} illustrates background theories related to developing the S-BORM algorithm. The algorithm is proposed in Section \ref{sec:sborm}. Performance of the algorithm is thoroughly investigated by three numerical examples in Section \ref{sec:ex} and various parametric studies in Section \ref{sec:param_study}. More technical details can be found in Appendix.

\section{Background}\label{sec:background}

\subsection{General system events and reliability-based optimization}

A failure event of a general system is represented as a series system of parallel systems of component failure events, or equivalently, a link-set of cut-sets, i.e.,
\begin{linenomath*}
\begin{equation}\label{eq:genSys_set}
    E_{\mathrm{sys}}=\bigcup_{k=1}^K \bigcap_{q\in\mathds{Q}_k} E_q,
\end{equation}
\end{linenomath*}
where $E_{\mathrm{sys}}$ and $E_q$ ($q\in\mathds{Q}_k$, $k =1,\ldots,K$) refer to the failure event of a system and component $q$, respectively. To determine whether a system or component event $E$ is either failure or survival, a limit-state function $g(\bm{x},\bm{V})$ can be used to account for performance of the corresponding system or component. The function depends on design variables $\bm{x}=(x_1,\ldots,x_D)$ and random vector $\bm{V}=(V_1,\ldots,V_M)$.  A realization $\bm{v}$ of $\bm{V}$ is considered a failure if $g(\bm{x},\bm{v})>0$ and a survival, otherwise\footnote{If necessary, one can reverse the definition of failure and survival (i.e., a failure event if $g(\bm{x},\bm{v})<0$) by reversing the sign of limit-state functions.}. Thereby, the definition of a system failure event in (\ref{eq:genSys_set}) can be represented in terms of limit-state functions as
\begin{linenomath*}
\begin{equation}\label{eq:genSys_lf}
    g_{\mathrm{sys}}(\bm{x},\bm{v}) = \max_{k={1,\ldots,K}} \min_{q\in\mathds{Q}_k} g_q(\bm{x},\bm{v}),
\end{equation}
\end{linenomath*}
where $g_{\mathrm{sys}}(\bm{x},\bm{v})$ and $g_q(\bm{x},\bm{v})$ denote the limit-state function of a system and a component, respectively. 

The most common formulation of reliability-based optimization is to minimize design cost $c(\bm{x})$ while satisfying a reliability constraint (i.e., system failure probability be less than a target value $p_f^t$):
\begin{linenomath*}
\begin{subequations} \label{rbo:basic}
\postdisplaypenalty=0
\begin{alignat}{2}
&\!\min_{\bm{x}\in\mathds{X}}        &\qquad& c(\bm{x})\\
&\text{subject to} &      & p(\bm{x}) \leq p_f^t,\label{rbo:basic_const}
\end{alignat}
\end{subequations}
\end{linenomath*}
where $\mathds{X}$ denotes a constraint set, and the conventional failure probability is defined as
\begin{linenomath*}
\begin{equation}\label{eq:pf_def}
    p(\bm{x})=P[g_{\mathrm{sys}}(\bm{x},\bm{V})>0].
\end{equation}
\end{linenomath*}
Alternatively, the reliability constraint (\ref{rbo:basic_const}) can be represented in terms of a quantile: $p(\bm{x}) \leq p_f^t$ if and only if $q_{1-p_f^t}(\bm{x}) \leq 0$, where $q_\alpha(\bm{x})$ is the $\alpha$-quantile of $g_{\mathrm{sys}}(\bm{x},\bm{V})$. Thus, (\ref{rbo:basic}) is equivalent to the problem
\begin{linenomath*}
\begin{subequations} \label{rbo:basic_q}
\postdisplaypenalty=0
\begin{alignat}{2}
&\!\min_{\bm{x}\in\mathds{X}}        &\qquad& c(\bm{x})\\
&\text{subject to} &      & q_{1-p_f^t}(\bm{x}) \leq 0.\label{rbo:basic_q_const}
\end{alignat}
\end{subequations}
\end{linenomath*}

It is noted that solving (\ref{rbo:basic}) and (\ref{rbo:basic_q}) requires an accessible expression for the reliability constraints (\ref{rbo:basic_const}) and (\ref{rbo:basic_q_const}). However, for general systems, deriving such an analytical formula is usually impossible. In such case, the failure probability in (\ref{eq:pf_def}) can be estimated by realizations of $\bm{V}$ (i.e., samples or data points), $\bm{v}_1,\ldots,\bm{v}_N$, leading to the formula
\begin{linenomath*}
\begin{equation} \label{eq:pf_estimate}
    \hat{p}(\bm{x}) = \sum_{n=1}^N p_n \cdot \mathds{I}[g_{\mathrm{sys}}(\bm{x},\bm{v}_n)>0],
\end{equation}
\end{linenomath*}
where $\hat{p}(\bm{x})$ is the estimated failure probability given $\bm{x}$, and $p_n$ is a probability or weight of realization $\bm{v}_n$\footnote{For example, if $\bm{v}_1,\ldots,\bm{v}_N$ are generated by Monte Carlo simulation (MCS), $p_n=1/N$ for all $n$.}. The Heaviside function $\mathds{I}[\cdot]$ takes value 1 if the given statement is true and 0, otherwise.
When replacing in \eqref{rbo:basic} the probability $p(\bm{x})$  by $\hat p(\bm{x})$, however, the Heaviside function in (\ref{eq:pf_estimate}) greatly complicates computation because it lacks (sub)gradient information. The gradient is not defined for $\bm{x}$ when there is some $n$ leading to $g_{\mathrm{sys}}(\bm{x},\bm{v}_n)=0$ and remains 0 at all other values.

\subsection{Buffered optimization and reliability method}\label{sec:borm}

The paper \cite{borm22} recently proposed a framework for reliability-based optimization, namely {\it buffered optimization and reliability method} (BORM). Being traced back to \cite{bpf10}, BORM replaces the failure probability in (\ref{rbo:basic_const}) by the buffered failure probability. While details can be found in the references, this section presents a brief illustration that is directly related to the following discussions.

The two failure probabilities are distinguished by how the threshold value of limit-state functions is defined to determine a failure event. According to the conventional probability, this threshold is fixed at 0. In contrast, the buffered probability defines such threshold as a quantile value whose associated superquantile is 0. In more detail, consider a random variable $Y$ and its CDF $F_Y(y)$. Then, the $\alpha$-quantile of $Y$, $q_\alpha$ is defined as
\begin{linenomath*}
\[
    q_\alpha=F_Y^{-1}(\alpha)
\]
\end{linenomath*}
provided that $F_Y$ is strictly increasing; a similar definition holds in general. The $\alpha$-superquantile, denoted by $\bar{q}_\alpha$, is defined as the average value of $Y$ beyond $q_\alpha$, i.e.,\footnote{If $Y$ is continuously distributed, then $\bar{q}_\alpha$ is equivalent to the conditional mean $\mathds{E}[Y|Y\geq q_\alpha]$.}
\begin{linenomath*}
\[
    \bar{q}_\alpha = q_\alpha + \frac{1}{1-\alpha}\mathds{E}[\max\{Y-q_\alpha,0\}].
\]
\end{linenomath*}
Finally, the buffered failure probability $\bar{p}_f$ of the event $Y>0$ is defined as
\begin{linenomath*}
\[
    \bar{p}_f = 1 - \bar{\alpha}_0,
\]
\end{linenomath*}
where $\bar{\alpha}_0$ is the probability that gives a zero superquantile, i.e., $\bar{q}_{\bar{\alpha}_0} = 0$.

These definitions motivate the shift from the RBO problem \eqref{rbo:basic} to the problem 
\begin{linenomath*}
\postdisplaypenalty=0
\begin{subequations}\label{rbo:basicB}
\begin{alignat}{2}
&\!\min_{\bm{x}\in\mathds{X}}        &\qquad& c(\bm{x})\\
&\text{subject to} &      & \bar p(\bm{x}) \leq \bar{p}_f^t,
\end{alignat}
\end{subequations}
\end{linenomath*}
where $\bar p(\bm{x})$ is the buffered failure probability of the event $g_{\mathrm{sys}}(\bm{x},\bm{V}) > 0$ and $\bar p_f^t$ is a threshold. We note that $\bar p(\bm{x}) \geq p(\bm{x})$ so that the buffered failure probability always bounds the conventional failure probability conservatively \cite{bpf10}. The new formulation facilitates a data-driven setting with $\bm{V}$ replaced by the outcomes $\bm{v}_1,\ldots,\bm{v}_N$. In this case, (\ref{rbo:basicB}) can be reformulated as
\begin{linenomath*}
\begin{subequations} \label{rbo:bpf}
\postdisplaypenalty=0
\begin{alignat}{2}
&\!\min_{\bm{x}\in\mathds{X},\gamma\in\mathds{R}}        &\qquad& c(\bm{x})\\
&\text{subject to} &      & \gamma + \frac{1}{\bar{p}_f^t} \sum_{n=1}^N p_n \max\{0,g_{\mathrm{sys}}(\bm{x},\bm{v}_n) - \gamma\} \leq 0\label{rbo:rel_const},
\end{alignat}
\end{subequations}
\end{linenomath*}
where $\gamma$ is an additional real-valued design variable, which at optimality specifies the $(1-\bar{p}_f^t)$-quantile value of $g_{\mathrm{sys}}(\bm{x},\bm{V})$. It is noted that since $\gamma$ appears in a well-structured manner, it does not increase the computational complexity of the optimization problem. The complexity is also not affected by the maximum operation in the constraint (\ref{rbo:rel_const}) as it can be reformulated to retain convexity and/or smoothness of the functions within the curly brackets (detailed illustrations can be found in \cite{bpf10}). Accordingly, optimization complexity is governed by $\mathds{X}$, $c(\bm{x})$, $g_\mathrm{sys}(\bm{x},\bm{v}_n)$, and sample size $N$. For instance, if these functions are convex and $\mathds{X}$ is a convex set, the problem becomes convex and is thus easily solvable using standard algorithms should $N$ be of moderate size. In such a convex setting, if instead $N$ is too large (say $N\geq 10^4$), then the problem can be efficiently solved by nonlinearly-constrained convex bundle methods such as the one proposed in \cite{vanAckooij_Oliveira_2014}. Even when the functions are neither linear nor convex, the gradients of component limit-state functions $g_q(\bm{x},\bm{v}_n)$ can be used for optimization algorithms, which greatly facilitates implementation as we see below.

In contrast to \cite{borm22}, which considers a ``system'' consisting of a single component, we now consider general systems.

\section{Proposed reliability-based optimization of general system events: S-BORM algorithm}\label{sec:sborm}

\subsection{Key ideas}

To develop an efficient optimization scheme that can handle general systems, we introduce four ideas for solving the data-driven RBO problem \eqref{rbo:bpf} in settings of general systems. First, the reliability constraint \eqref{rbo:rel_const} is penalized and moved to the objective function. Second, we linearlize the limit-state functions adaptively at the current candidate solution. Third, we reformulate the now modified objective as a difference-of-convex function, which produces a subproblem solvable by the difference-of-convex bundle method of \cite{pb17}. Fourth, we improve computational efficiency further by employing an active-set strategy.
That is, at each iteration of optimization, the algorithm considers only a subset of samples that are within or close enough to failure domains at a current solution, i.e., the samples with the highest limit-state function values. Since the number of failure events are in general very small, this strategy greatly facilitates optimization. More details about active-set strategies are available in \cite{borm22}. We name the proposed approach the S-BORM algorithm as it is designed to handle system events. 

Our approach follows the standard path consisting of replacing a difficult problem with a sequence of simpler subproblems and/or models. However, our pathway is distinguished from what is largely considered in the optimization literature, where the subproblems typically have easily obtainable solutions via quadratic or convex optimization. We utilize more complex subproblems, which appears necessary to capture the max-min formula \eqref{eq:genSys_lf} for the system limit-state function. In turn, this requires us to adopt more advanced subroutines for solving the subproblems. Specifically, we leverage Algorithm 1 in \cite{pb17}. The next subsections discuss the S-BORM algorithm in detail.

\subsection{Linearization of limit-state functions for difference-of-convex decomposition}\label{sec:lin_g}

We introduce mild conditions on RBO problems: $\mathds{X}$ is a polyhedral set and $c(\bm{x})$ as well as $g_q(\bm{x})$ are smooth (in $\bm{x}$) with Lipschitz continuous gradients on $\mathds{X}$. These conditions should hold for many practical problems. It is noted that convexity is not assumed for either cost function or limit-state functions.

By recalling (\ref{eq:genSys_lf}), the optimization problem (\ref{rbo:bpf}) takes the form
\begin{linenomath*}
\begin{subequations} \label{rbo:bpf_genSys}
\postdisplaypenalty=0
\begin{alignat}{2}
&\!\min_{\bm{x}\in\mathds{X},\gamma\in \mathds{R}}        &\qquad& c(\bm{x})\\
&\text{subject to} &      & \gamma + \frac{1}{\bar p_f^t} \sum_{n=1}^N p_n \max\{0, \max_{k\in{1,\ldots,K}} \min_{q\in\mathds{Q}_k} g_q(\bm{x},\bm{v}_n)-\gamma \} \leq 0.
\end{alignat}
\end{subequations}
\end{linenomath*}
We penalize the reliability constraint: for $\theta\in(0,\infty)$, the optimization problem becomes 
\begin{linenomath*}
\begin{equation} \label{rbo:bpf_penalize}
\min_{\bm{x}\in\mathds{X},\gamma\in \mathds{R}} 
F(\bm{x},\gamma;\theta),\; \mbox{ with }\; F(\bm{x},\gamma;\theta):=c(\bm{x}) + \theta \max\bigg\{ 0, \gamma + \frac{1}{\bar{p}_f^t} \sum_{n\in\hat{\mathds{N}}} p_n \max\Big\{0, \max_{k=1,\ldots,K} \min_{q\in\mathds{Q}_k} g_q(\bm{x},\bm{v}_n)-\gamma \Big\} \bigg\}.
\end{equation}
\end{linenomath*}
Above, an active-set strategy is assumed, i.e., the problem considers only active samples $\bm{v}_n$ with $n\in\hat{\mathds{N}}\subset\{1,\ldots,N\}$. The problem is further approximated by linearizing each $g_q(\bm{x},\bm{v}_n)$ at a candidate solution $\hat{\bm{x}}^\nu$:
\begin{linenomath*}
\begin{equation} \label{rbo:bpf_penal_linear}
    \min_{\bm{x}\in\mathds{X},\gamma\in \mathds{R}} \quad F^\nu(\bm{x},\gamma;\theta,\hat{\bm{x}}^\nu),
\end{equation}
\end{linenomath*}
where
\begin{linenomath*}
\begin{equation}\label{eq:obj_pen}
    F^\nu(\bm{x},\gamma;\theta,\hat{\bm{x}}^\nu) := c(\bm{x}) + \theta \max\bigg\{ 0, \gamma + \frac{1}{\bar{p}_f^t} \sum_{n\in\hat{\mathds{N}}} p_n \max\Big\{0, \max_{k\in{1,\ldots,K}} \min_{q\in\mathds{Q}_k} g_q(\hat{\bm{x}}^\nu,\bm{v}_n) + \langle \nabla g_q(\hat{\bm{x}}^\nu,\bm{v}_n), \bm{x}-\hat{\bm{x}}^\nu \rangle -\gamma \Big\} \bigg\}.
\end{equation}
\end{linenomath*}

The function $F^\nu$ can be written as a difference-of-convex function, i.e., a convex function minus another convex function. We derive the specific formula in \ref{sec:dc_deriv}. Thus, \eqref{rbo:bpf_penal_linear} is a subproblem that can be addressed by Algorithm 1 in \cite{pb17}. The S-BORM algorithm solves such subproblems, with slight adjustments, repeatedly as described next.

\subsection{S-BORM algorithm}\label{sec:sborm_algorithm}

Based on the derivation in Sections \ref{sec:lin_g} and \ref{sec:dc_deriv}, we now present the S-BORM algorithm.

{\bf S-BORM Algorithm:}
\begin{itemize}[itemsep=0pt, topsep=0pt, label=\textbf{Data}, leftmargin=1.5cm]
    \item Given an initial point $\bm{x}^0\in\mathds{X}$, samples $\{\bm{v}_1,\ldots,\bm{v}_N\}$, and a parameter $\gamma^0$, choose algorithm parameters $\theta>0$, $\theta^{\max}>\theta$, $\lambda>0$, $\omega\geq 1$, $\kappa\in(0,1)$, and {\tt tol}.
\end{itemize}
\begin{enumerate}[itemsep=0pt, topsep=0pt, label=\textbf{Step \arabic*}, leftmargin=1.75cm, start=0]
    \item Set $\nu=0$ and $\hat{\bm{x}}^\nu=\bm{x}^0$, $\hat{\gamma}^\nu=\gamma^0$, $\theta^\nu = \theta$, $\lambda^\nu = \lambda$.
    
    \item Evaluate $g_{\mathrm{sys}}(\hat{\bm{x}}^\nu,\bm{v}_1),\ldots,g_{\mathrm{sys}}(\hat{\bm{x}}^\nu,\bm{v}_N)$ and obtain an index set of active samples, $\hat{\mathds{N}}^\nu \subset \{1,\ldots,N\}$, with the $\lceil \omega N \bar{p}_f^t \rceil$ greatest values of $g_{\mathrm{sys}}(\hat{\bm{x}}^\nu,\bm{v}_n)$, $n=1,\dots,N$.\label{sborm_al:g_eval}
    
    \item Compute $\nabla g_q(\hat{\bm{x}}^\nu,\bm{v}_n)$ and $g_q(\hat{\bm{x}}^\nu,\bm{v}_n)$ for all $q\in\mathds{Q}_k$, $k=1,\ldots,K$ and $n \in \hat{\mathds{N}}^\nu$.\label{sborm_al:linearize}
    
    \item\label{sborm_al:pb_eval} 
    Let $F^\nu$ as in \eqref{eq:obj_pen} and solve the subproblem
    \begin{linenomath*}
    \begin{equation}\label{rbo:in_sborm}
        \min_{\bm{x}\in \mathds{X}, \gamma\in\mathds{R}} F^\nu(\bm{x},\gamma;\theta^\nu,\hat{\bm{x}}^\nu) + \frac{\lambda^\nu}{2} \|\bm{x}-\hat{\bm{x}}^\nu\|^2_2 + \frac{\lambda^\nu}{2} (\gamma - \hat{\gamma}^\nu)^2
    \end{equation}
    \end{linenomath*}
    and obtain a point $(\bm{x}^{\nu+1},\gamma^{\nu+1})$ using \cite[Alg. 1]{pb17} with initial point $(\hat{\bm{x}}^{\nu},\hat{\gamma}^{\nu})$.
    If $\norm{\bm{x}^{\nu+1}-\hat{\bm{x}}^\nu}^2_2 + (\gamma^{\nu+1} - \hat{\gamma}^\nu)^2 \leq \tt{tol}$, stop and return $(\hat{\bm{x}}^{\nu},\hat{\gamma}^{\nu})$ as a potential solution. Else, go to \ref{sborm_al:eval_sol}.
    
    \item\label{sborm_al:eval_sol} With $F$ from \eqref{rbo:bpf_penalize}, define the predicted decrease as
    \begin{linenomath*}
    \begin{equation*}
        \zeta^\nu := F(\hat{\bm{x}}^{\nu},\hat{\gamma}^{\nu};\theta^\nu)-F^\nu(\bm{x}^{\nu+1},\gamma^{\nu+1};\theta^\nu,\hat{\bm{x}}^\nu) - \frac{\lambda^\nu}{2} \|\bm{x}^{\nu+1}-\hat{\bm{x}}^\nu\|^2_2 - \frac{\lambda^\nu}{2} (\gamma^{\nu+1}-\hat{\gamma}^\nu)^2 .
    \end{equation*}
    \end{linenomath*}
    If $F(\bm{x}^{\nu+1},\gamma^{\nu+1};\theta^\nu) \leq F(\hat{\bm{x}}^{\nu},\hat{\gamma}^{\nu};\theta^\nu) - \kappa \zeta^\nu$, declare a \textbf{Serious Step:}
    \begin{linenomath*}
    \begin{equation*}
        (\hat{\bm{x}}^{\nu+1},\hat{\gamma}^{\nu+1}) := (\bm{x}^{\nu+1},\gamma^{\nu+1}) \quad \textrm{and} \quad \lambda^{\nu+1} := \lambda^\nu.
    \end{equation*}
    \end{linenomath*}
    Else, declare a \textbf{Null Step:}
    \begin{linenomath*}
    \begin{equation*}
        (\hat{\bm{x}}^{\nu+1},\hat{\gamma}^{\nu+1}) := (\hat{\bm{x}}^\nu,\hat{\gamma}^\nu) \quad \textrm{and} \quad \lambda^{\nu+1} := 2\lambda^\nu.
    \end{equation*}
    \end{linenomath*}
    
    \item\label{sborm_al:update_theta} Set $\theta^{\nu+1}:=\min \{ 1.5\theta^\nu, \theta^{\max} \}$. Replace $\nu$ by $\nu+1$ and go to \ref{sborm_al:g_eval} if Serious Step or \ref{sborm_al:pb_eval} if Null Step.
\end{enumerate}

\bigbreak 
The procedure of the algorithm is as follows. \ref{sborm_al:g_eval} evaluates limit-state functions with a current candidate solution $\hat{\bm{x}}^\nu$ and sorts out active samples $\hat{\mathds{N}}$. Since there are $N\bar{p}_f^t$ samples of system failure when a failure probability constraint is satisfied, the size of $\hat{\mathds{N}}$ is set as $\lceil \omega N \bar{p}_f^t \rceil$ with a parameter $\omega \geq 1$. Then, in \ref{sborm_al:linearize}, the limit-state functions in the penalized problem~\eqref{rbo:bpf_penalize} are linearized to define a subproblem. \ref{sborm_al:pb_eval} defines the next iterate as a solution of the subproblem\footnote{It suffices for the solution to be a critical point of the subproblem.}. If the new iterate is close enough to a current candidate solution, the algorithm is terminated. Otherwise, in \ref{sborm_al:eval_sol}, a descent test is carried out to ensure that the algorithm makes progress in each serious step. If the new iterates lead to satisfactory decrease in the objective function, this update is considered a Serious Step: the new iterate becomes the current candidate solution. Otherwise, it is a Null Step: the candidate solution does not change, and the parameter $\lambda^\nu$ is increased to force the next iterate to be closer to the previous point $\hat x^\nu$ so far. At the end of each iteration, in \ref{sborm_al:update_theta}, $\theta^\nu$ that penalizes violating reliability constraints is increased. Thereby, solutions can be gradually nudged towards sufficient reliability.

The proposed algorithm has three major parameters: the two penalty terms $\lambda$ and $\theta$ and the ratio of active samples, $\omega$. The default values of these parameters are proposed as $\lambda=0.01$, $\theta=1$, $\theta^{\max}=10^5$ and $\omega=2$, which are also used in the following examples. Other minor parameters are $\kappa$ and {\tt tol}; $\kappa$ decides the criterion whether to accept a new candidate solution, and {\tt tol} is the stopping threshold of Euclidean distance between the candidate solution and the new iterate. Default values of both parameters are proposed as 0.01. It is noted that these parameters have little influence on optimization results as demonstrated by numerical examples in Section \ref{sec:param_study}, where we show robustness of the proposed algorithm by performing parametric study under various settings. 

\section{Numerical examples}\label{sec:ex}

\subsection{Experiment settings}\label{sec:ex_param}

We design three numerical examples to demonstrate how the S-BORM algorithm provides an efficient and accessible means to solve problems that have been considered challenging. The first two examples are structural systems. A common way to define failure of a structural system is to identify multiple failure modes. Then, a system failure is defined as occurrence of any of the failure modes. However, because of computational difficulty, failure modes are often handled separately by each being assigned a target failure probability \cite{rbo-kriging21,borm22}. In contrast, the S-BORM algorithm handles system failure without simplification. In addition, the second example shows the capability of the algorithm for handling a large number of failure modes.

The third example investigates the issue of identifying optimal allocations of testing time across components from the perspective of system reliability. This requires us to combine reliability optimization with reliability growth models (RGMs). Although the issue has been investigated from various aspects (e.g., identification of most critical components \cite{rgm10} or adaptive testing strategies \cite{rgm19}), this example is the first attempt to perform optimization within the context of a general system.

Algorithm parameters are set as the default values proposed in Section \ref{sec:sborm_algorithm}. For all problems, initial solutions are set as the midpoint of upper and lower bounds of design variables. The target buffered failure probability is $\bar{p}_f^t=1 \cdot 10^{-3}$, and the target coefficient of variance (c.o.v.) is $\delta^t=0.05$. This leads to the number of samples, $(1-\bar{p}_f^t)/(\bar{p}_f^t \cdot (\delta^t)^2) = 399{,}600$ \cite{borm22}. Then, with the default parameter $\omega=2$, the number of active samples becomes $\lceil 2 \cdot 399{,}600 \cdot 10^{-3} \rceil=800$.

\subsection{Design of cantilever beam-bar system}

\begin{figure}[h!]
    \centering
    \includegraphics[scale = 0.5]{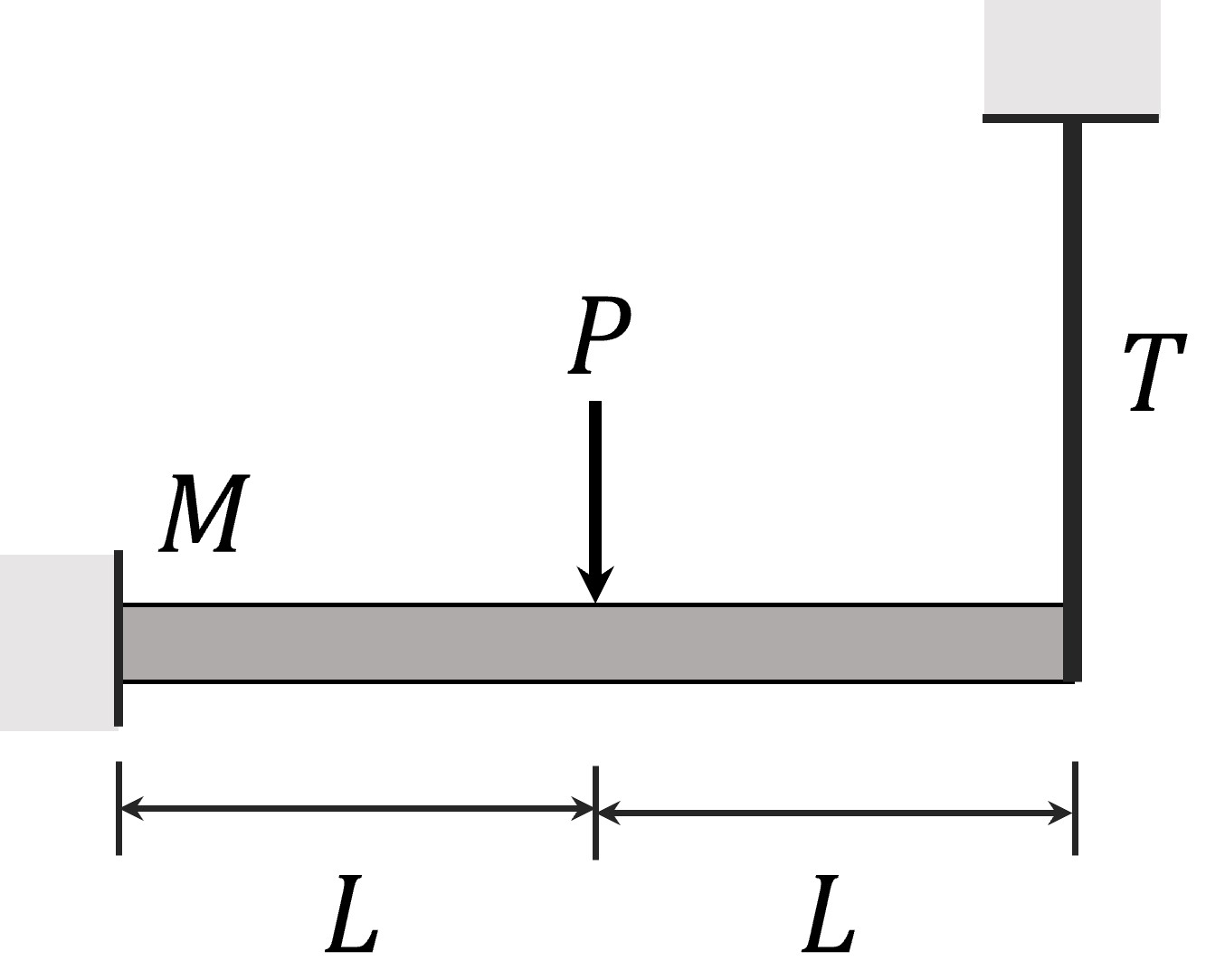}
    \caption{Example cantilever beam-bar system (figure recreated from \cite{LpBound03})}
    \label{fig:beam}
\end{figure}

This example investigates an optimal design of a cantilever beam-bar system illustrated in Figure \ref{fig:beam}. The system consists of an ideally plastic cantilever beam of moment capacity $M$ and length $2L=2\cdot5$, which is propped by an ideally rigid-brittle bar of strength $T$ \cite{LpBound03}. The structure is subjected to load $P$ that is applied on the middle of the bar. Load $P$ is a random variable following the normal distribution with mean $\mu_P=150$ and standard deviation $\sigma_P=30$. Moment $M$ and strength $T$ are also normal random variables with mean $\mu_M$ and $\mu_T$ and standard deviation $\sigma_M=300$ and $\sigma_T=20$, respectively. In this example, we optimize two design variables $x_1 = \mu_M \in [500, 1500]$ and $x_2 = \mu_T\in[50, 150]$. The cost function is $c(x)=2x_1+x_2$.

As seen from \cite{LpBound03}, we identify six limit-state functions such that 
\begin{linenomath*}
\postdisplaypenalty=0
\begin{align*}
    & g_1( \bm{x}, \bm{v} ) = -(x_2 + v_2 - 5 v_3/16), \\
    & g_2( \bm{x}, \bm{v} ) = -(x_1 + v_1 - L v_3), \\
    & g_3( \bm{x}, \bm{v} ) = -(x_1 + v_1 - 3L v_3/8), \\
    & g_4( \bm{x}, \bm{v} ) = -(x_1 + v_1 - L v_3/3), \\
    & g_5( \bm{x}, \bm{v} ) = -(x_1 + v_1 + 2L \cdot (x_2 + v_2) - L v_3),
\end{align*}
\end{linenomath*}
where $v_1$ and $v_2$ are realizations of the normal distribution with zero mean and standard deviation $\sigma_M$ and $\sigma_T$, respectively; and $v_3$ is a realization of $P$. The system event consists of three cut-sets such that $E_{\mathrm{system}}=E_1 E_2 \cup E_3 E_4 \cup E_3 E_5$. 

Using the S-BORM algorithm, the computed solution is $(x_1,x_2)=(1297, 150.0)$, resulting in estimated buffered failure probability ${\hat {\bar p}}_f=9.985 \cdot 10^{-4}$ and cost 2,743. For comparison, a grid search is performed by discretizing each variable into 100 intervals, which leads to 10.1 and 1.01 spacing for $x_1$ and $x_2$, respectively. The search identifies the best point as (1298, 149.0) leading to cost 2,745, ${\hat{\bar p}}_f=9.935\cdot 10^{-4}$, and ${\hat p}_f=2.678\cdot 10^{-4}$. This agrees with the solution computed by the proposed algorithm. The optimization requires a marginal computational efforts taking 0.1 seconds. It runs 7 outer loops and 7 rounds of gradient evaluations, which implies that only 7 evaluations of limit-state functions are made for each sample and 7 gradient evaluations are made for each active sample. A summary of the results can be found in Table \ref{tab:ex_results}.

\subsection{Design of indeterminate truss bridge system}

\begin{figure}[h!]
    \centering
    \includegraphics[scale = 0.6]{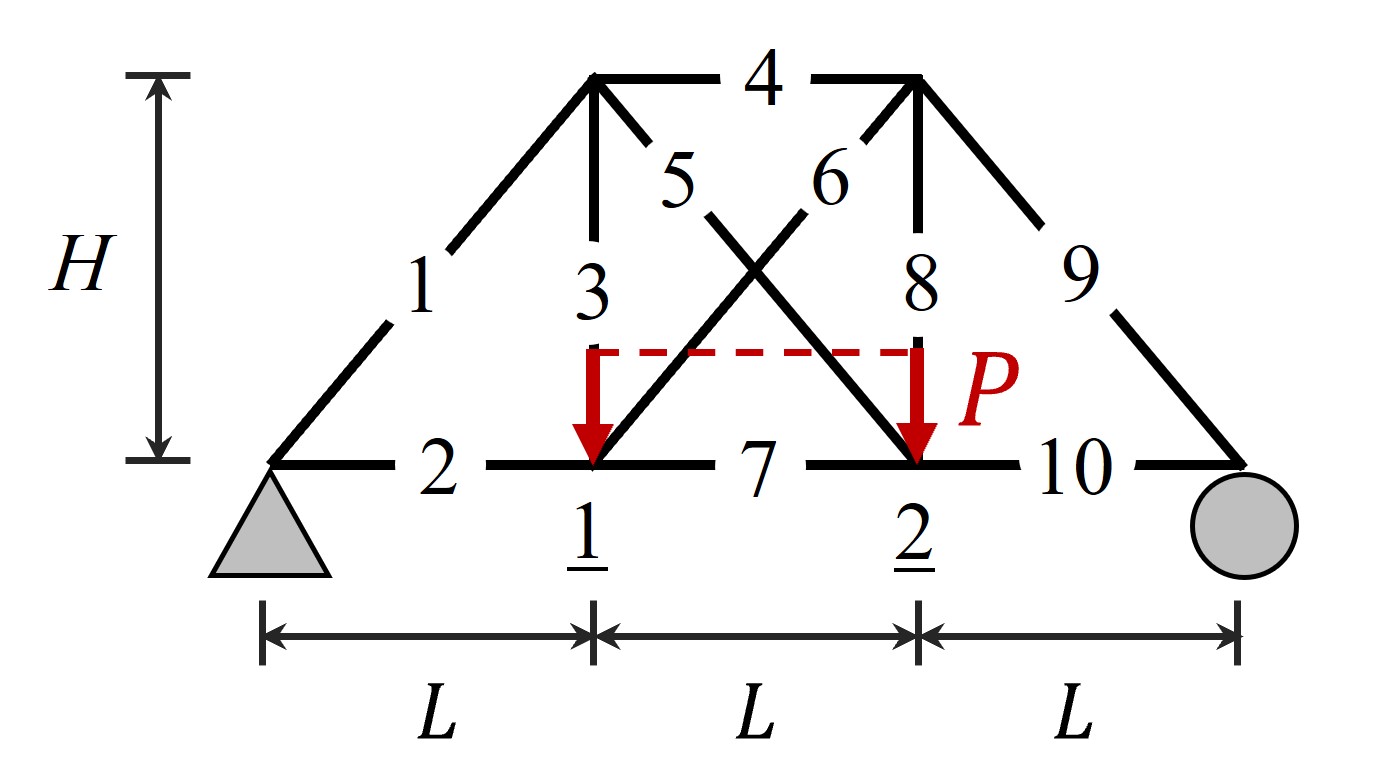}
    \caption{Example truss structure system}
    \label{fig:truss}
\end{figure}

Consider the truss bridge structure in Figure \ref{fig:truss}. The structure consists of 10 members with strength $R_i$, $i=1,\ldots,10$, which are independent normal random variables with mean $\mu_R = 276$ and standard deviation $\sigma_R = 13.8$. The parameters of member lengths are set as $H=1.6$ and $L=2$. The truss structure is subjected to load $P$ exerted on nodes 1 and 2, which follows the normal distribution with mean $\mu_P=190$ and standard deviation $\sigma_P = 19$. Then, the design variables are cross-section areas of the members. To reflect the practical aspect of construction, the members are grouped into four sets such that $\{1,2,9,10\}$, $\{3,8\}$, $\{4,7\}$, and $\{5,6\}$. Within a group, members have the same cross-section area, which is represented by a design variable $x_d \in [1,2]$, $d=1,\ldots,4$. The cost function is the total volume of the members, i.e., $c(\bm{x})=\sum_{q=1}^{10} l_q \cdot x_{d(q)}$, where $l_q$ is the length of member $q$, and $d(q)$ denotes the index of the group that member $q$ belongs to, e.g., $d(1)=1$ and $d(3)=2$.

The system failure is defined as an occurrence of structural instability. Such occurrence can be identified by performing structural analysis. Each failure mode, i.e., a sequence of member failures that leads to a system failure, can be represented as a cut-set, whereby the system event is defined as a link-set of those cut-sets. For example, as a failure of member 1 incurs structural instability, member 1 constitutes a cut-set. Another example is a failure of member 3 followed by that of member 7. Such failure of member $q$ at failure mode $k$ can be represented by a limit-state function
\begin{linenomath*}
\begin{equation*}
    g_{kq}(\bm{x},\bm{v}) = v_0 \cdot \delta_{kq} - x_{d(q)} \cdot v_q,
\end{equation*}
\end{linenomath*}
where $v_0$ is a realization of $P$; $\delta_{kq}$ is the force experienced by member $q$ in failure mode $k$ when a unit force is applied to nodes 1 and 2; and $v_q$, $q=1,\ldots,10$, is a realization of $R_q$. As 108 failure modes are identified for this structure, the system event becomes
\begin{linenomath*}
\begin{equation*}
    E_{\mathrm{sys}} = \bigcup_{k=1}^{108} \bigcap_{q \in \mathds{Q}_k} E_{kq},
\end{equation*}
\end{linenomath*}
where $\mathds{Q}_k$ denotes the index set of members that constitute a failure mode $k$.

As summarized in Table \ref{tab:ex_results}, a solution is obtained as (1.586, 1.000, 1.459, 1.000) with cost 28.63, ${\hat{\bar p}}_f=9.735\cdot 10^{-4}$, and ${\hat p}_f=3.654\cdot 10^{-4}$. The optimization takes 13.6 seconds with 3 rounds of outer loops and 3 rounds of gradient evaluations.

\subsection{Testing time allocation on electrical components}

\begin{figure}[h!]
    \centering
    \includegraphics[scale = 0.5]{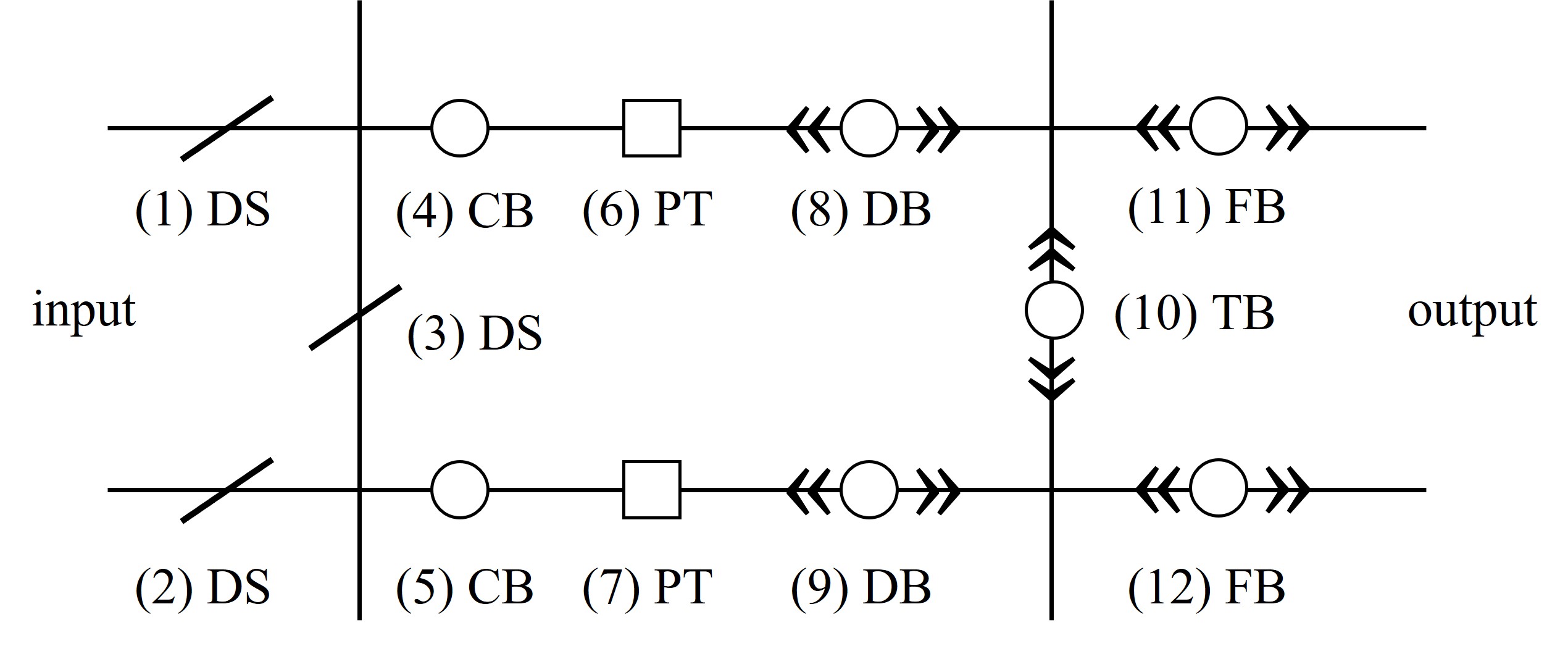}
    \caption{Example electrical system (figure recreated from \cite{power08})}
    \label{fig:power}
\end{figure}

Consider the two transmission-line electrical substation system in Figure \ref{fig:power} \cite{power08}. The system consists of 12 components, which either fail or survive. Then, system failure is defined as a disconnection between the input and output nodes. Each component type is under a test phase, and we aim to find an optimal allocation of testing time over the component types. There are 6 component types, i.e., disconnect switch (DS), circuit breaker (CB), power transformer (PT), drawout breaker (DB), tie breaker (TB), and feeder breaker (FB), and their testing time is denoted by design variables $x_1,\ldots,x_6$, respectively. We assume that their reliability growth follows the non-homogeneous Poisson process (NHPP) RGM \cite{MsrRgm17}, by which the fault rate of component $q$ of type $d(q)$, $q=1,\ldots,12$ is defined as
\begin{linenomath*}
\begin{equation*}
    \lambda_q(x_{d(q)}) = \frac{\alpha \beta}{ \mathrm{exp}( \beta x_{d(q)} ) },
\end{equation*}
\end{linenomath*}
where the parameters are set as $\alpha=9$ and $\beta=2$. The cost function is the total testing time, i.e., $c(x)=\sum_{d=1}^{4} x_d$. On the other hand, for each component $q$, the limit-state function is defined as
\begin{linenomath*}
\postdisplaypenalty=0
\begin{align*}
    g_q(x_{d(q)},v_q) & = {\Delta T - \frac{-\mathrm{ln} v_q}{\lambda_q (x_{d(q)})}}  \\
    & = {\Delta T + \mathrm{ln} v_q \cdot \frac{\mathrm{exp}(\beta x_{d(q)})} { \alpha \beta  }},
\end{align*}
\end{linenomath*}
where $\Delta T=365$ is the target operation time, and $\xi_q$, $q=1,\ldots,12$, is a realization of the uniform distribution $U[0,1]$. The system consists of 25 minimum cut sets \{(1, 2), (4, 5), (4, 7), (4, 9), (5, 6), (6, 7), (6, 9), (5, 8), (7, 8), (8, 9), (11, 12), (1, 3, 5), (1, 3, 7), (1, 3, 9), (2, 3, 4), (2, 3, 6), (2, 3, 8), (4, 10, 12), (6, 10, 12), (8, 10, 12), (5, 10, 11), (7, 10, 11), (9, 10, 11), (1, 3, 10, 12), (2, 3, 10, 11)\} \cite{power08}.

As summarized in Table \ref{tab:ex_results}, a solution is computed as (7.017,  7.047, 7.095, 7.024, 1.000, 7.016) with cost 36.20, ${\hat{\bar p}}_f=9.635\cdot 10^{-4}$, and ${\hat p}_f=4.179\cdot 10^{-4}$. Since all component types are assigned an identical RGM, difference in testing time arises solely from varying topological importance. The optimization takes 2.83 seconds with 10 rounds of evaluations of limit-state functions and 10 rounds of gradient evaluations of active samples.

\begin{table}[h!]
    \centering
    \begin{tabular}{c|c|c|c|c|c|c|c}
         \hline
         Ex. & \makecell{Time \\[-5pt] (sec.)} & \makecell{No. of outer loops \\[-5pt] (function calls)$^\mathrm{a}$} & \makecell{No. of serious steps \\[-5pt] (gradient calls)$^\mathrm{b}$} & \makecell{Computed \\[-5pt] solution} & Cost & \makecell{${\hat {\bar p}}_f$ \\[-5pt] $(\cdot 10^{-4})$} & \makecell{${\hat p}_f$ \\[-5pt] $(\cdot 10^{-4})$} \\
         \hline
         1 & 0.136 & \makecell{7 \\[-5pt] (2,797,200)} & \makecell{7 \\[-5pt] (5,600)} & (1297, 150.0) & 2,743 & 9.985 & 2.678 \\
         \hline
         2 & 13.6 & \makecell{3 \\[-5pt] (1,198,800)} & \makecell{3 \\[-5pt] (2,400)} & \makecell{(1.586, 1.000, \\[-5pt] 1.459, 1.000)} & 28.63 & 9.735 & 3.654 \\
         \hline
         3 & 2.83 & \makecell{10 \\[-5pt] (3,996,000)} & \makecell{10 \\[-5pt] (8,000)} & \makecell{(7.017, 7.047, \\[-5pt] 7.095, 7.024, \\[-5pt] 1.000, 7.016)} & 36.20 & 9.860 & 4.429 \\
         \hline
         \multicolumn{8}{l}{\small $^\mathrm{a}$ (No. of outer loops)$\cdot$(No. of samples) $=$ (Total evaluation number of limit-state functions)} \\[-6pt]
         \multicolumn{8}{l}{\small $^\mathrm{b}$ (No. of serious steps)$\cdot$(No. of active samples) $=$ (Total evaluation number of  limit-state function gradients)}
    \end{tabular}
    \caption{Optimization results of three numerical examples}
    \label{tab:ex_results}
\end{table}

\section{Parametric test of the proposed algorithm}\label{sec:param_study}

\subsection{Algorithm parameters}

While default values of the algorithm parameters are proposed in Section \ref{sec:ex_param}, we test the robustness of the algorithm by running optimization with their values changed. First, experiments are performed by changing $\lambda$ to 0.005, 0.02, 0.04, 0.08, and 1. The results are summarized in Table \ref{tab:param_lam}. As illustrated in the table, the parameter does not incur notable differences. In all of the three examples, the computation time remains stable taking less than 1 minute. Also, for other results including the number of outer loops and gradient evaluations, cost, ${\hat{\bar p}}_f$, and $\hat{p}_f$, highest and lowest values all remain very close. Similarly, $\theta$ is changed to 0.25, 0.5, 2, 4, and 8. Table \ref{tab:param_theta} summarizes the highest and lowest values of various results. Again, the results do not show notable variances. Finally, $\omega$ is tested with values 1.2, 1.5, 2, 3, and 5, as illustrated in Table \ref{tab:param_omega}. This parameter also does not lead to meaningful differences.

Although the parameters are found to have insignificant influences, there are still marginal variances in computation time and quality of solutions (i.e., the cheapest solution that satisfies reliability constraints). Numerical experiments suggest that the proposed default values show the most stable performance. Nevertheless, one may alter their values to fit the characteristics of a given problem. For example, one may increase $\lambda$ if solutions are too slow to converge, jumping between distant regions. The parameter $\theta$ can be increased if the algorithm finds it hard to satisfy reliability constraints. The parameter $\omega$ needs to be increased if active sets change wildly as an optimal solution is updated, or decreased if evaluation of limit-state functions is costly.

\begin{table}[h!]
    \centering
    \begin{tabular}{c|c|c|c|c|c|c|c}
         \hline
         \multicolumn{2}{c|}{Ex.} & \makecell{Time \\[-5pt] (sec.)} & \makecell{No. of \\[-5pt] outer loops} & \makecell{No. of \\[-5pt]serious steps} & Cost & \makecell{${\hat {\bar p}}_f$ \\[-5pt] ($\cdot 10^{-4}$) } & \makecell{ ${\hat p}_f$ \\[-5pt] ($\cdot 10^{-4}$) } \\
         \hline
         \multirow{ 2}{*}{ \makecell{ \\[-5pt] 1 } } & Highest & \makecell{0.73 \\[-5pt] ($\lambda=1$)} & \makecell{26 \\[-5pt] ($\lambda=1$)} & \makecell{26 \\[-5pt] ($\lambda=1$)} & \makecell{2,743 \\[-5pt] ($\lambda=0.005$)} & \makecell{ 9.985 \\[-5pt] (all $\lambda$) } & \makecell{ 3.103 \\[-5pt] ($\lambda=1$) } \\ \cline{2-8}
         & Lowest & \makecell{0.14 \\[-5pt] ($\lambda=0.01$)} & \makecell{7 \\[-5pt] ($\lambda=0.01$)} & \makecell{7 \\[-5pt] ($\lambda=0.01$)} & \makecell{2,718 \\[-5pt] ($\lambda=0.04$)} & \makecell{ 9.985 \\[-5pt] (all $\lambda$) } & \makecell{ 2.678 \\[-5pt] ($\lambda=0.005$) } \\
         \hline
         \multirow{ 2}{*}{ \makecell{ \\[-5pt] 2 } } & Highest & \makecell{55 \\[-5pt] ($\lambda=1$)} & \makecell{4 \\[-5pt] ($\lambda=0.05$)} & \makecell{4 \\[-5pt] ($\lambda=0.05$)} & \makecell{29.35 \\[-5pt] ($\lambda=0.04$)} & \makecell{ 9.760 \\[-5pt] ($\lambda=1$) } & \makecell{ 3.754 \\[-5pt] ($\lambda=1$) } \\ \cline{2-8}
         & Lowest & \makecell{19 \\[-5pt] ($\lambda=0.01$)} & \makecell{2 \\[-5pt] ($\lambda=1$)} & \makecell{2 \\[-5pt] ($\lambda=1$)} & \makecell{28.63 \\[-5pt] ($\lambda=0.01$)} & \makecell{ 8.834 \\[-5pt] ($\lambda=0.08$) } & \makecell{ 3.353 \\[-5pt] ($\lambda=0.005$) } \\
         \hline
         \multirow{ 2}{*}{ \makecell{ \\[-5pt] 3 } } & Highest & \makecell{4.9 \\[-5pt] ($\lambda=1$)} & \makecell{10 \\[-5pt] (all $\lambda$)} & \makecell{10 \\[-5pt] (all $\lambda$)} & \makecell{36.78 \\[-5pt] ($\lambda=1$)} & \makecell{ 9.960 \\[-5pt] ($\lambda=0.08$) } & \makecell{ 4.354 \\[-5pt] ($\lambda=0.08$) } \\ \cline{2-8}
         & Lowest & \makecell{2.1 \\[-5pt] ($\lambda=0.08$)} & \makecell{10 \\[-5pt] (all $\lambda$)} & \makecell{10 \\[-5pt] (all $\lambda$)} & \makecell{36.14 \\[-5pt] ($\lambda=0.08$)} & \makecell{ 9.635 \\[-5pt] ($\lambda=0.04$) } & \makecell{ 4.304 \\[-5pt] ($\lambda=0.04$) } \\
         \hline
    \end{tabular}
    \caption{Optimization results obtained by different parameters $\lambda=0.005, 0.01, 0.02, 0.04, 0.08,$ and 1}
    \label{tab:param_lam}
\end{table}

\begin{table}[h!]
    \centering
    \begin{tabular}{c|c|c|c|c|c|c|c}
         \hline
         \multicolumn{2}{c|}{Ex.} & \makecell{Time \\[-5pt] (sec.)} & \makecell{No. of \\[-5pt] outer loops} & \makecell{No. of \\[-5pt]serious steps} & Cost & \makecell{${\hat {\bar p}}_f$ \\[-5pt] ($\cdot 10^{-4}$) } & \makecell{ ${\hat p}_f$ \\[-5pt] ($\cdot 10^{-4}$) } \\
         \hline
         \multirow{ 2}{*}{ \makecell{ \\[-5pt] 1 } } & Highest & \makecell{0.21 \\[-5pt] ($\theta=0.05$)} & \makecell{12 \\[-5pt] ($\theta=0.25$)} & \makecell{10 \\[-5pt] ($\theta=0.5$)} & \makecell{2,743 \\[-5pt] ($\theta=0.25$)} & \makecell{ 9.985 \\[-5pt] (all $\theta$) } & \makecell{ 3.128 \\[-5pt] ($\theta=8$) } \\ \cline{2-8}
         & Lowest & \makecell{0.14 \\[-5pt] ($\theta=1$)} & \makecell{3 \\[-5pt] ($\theta=4$)} & \makecell{3 \\[-5pt] ($\theta=4$)} & \makecell{2,718 \\[-5pt] ($\theta=4$)} & \makecell{ 9.985 \\[-5pt] (all $\theta$) } & \makecell{ 2.678 \\[-5pt] ($\theta=0.25$) } \\
         \hline
         \multirow{ 2}{*}{ \makecell{ \\[-5pt] 2 } } & Highest & \makecell{51 \\[-5pt] ($\theta=0.5$)} & \makecell{4 \\[-5pt] ($\theta=0.25$)} & \makecell{3 \\[-5pt] ($\theta=8$)} & \makecell{28.82 \\[-5pt] ($\theta=2$)} & \makecell{ 9.960 \\[-5pt] ($\theta=2$) } & \makecell{ 3.879 \\[-5pt] ($\theta=2$) } \\ \cline{2-8}
         & Lowest & \makecell{19 \\[-5pt] ($\theta=1$)} & \makecell{2 \\[-5pt] ($\theta=2$)} & \makecell{2 \\[-5pt] ($\theta=2$)} & \makecell{28.61 \\[-5pt] ($\theta=8$)} & \makecell{ 7.432 \\[-5pt] ($\theta=4$) } & \makecell{ 2.928 \\[-5pt] ($\theta=4$) } \\
         \hline
         \multirow{ 2}{*}{ \makecell{ \\[-5pt] 3 } } & Highest & \makecell{2.2 \\[-5pt] ($\theta=2$)} & \makecell{10 \\[-5pt] (all $\theta$)} & \makecell{10 \\[-5pt] (all $\theta$)} & \makecell{39.21 \\[-5pt] ($\theta=8$)} & \makecell{ 9.910 \\[-5pt] ($\theta=2$) } & \makecell{ 4.455 \\[-5pt] ($\theta=2$) } \\ \cline{2-8}
         & Lowest & \makecell{2.1 \\[-5pt] ($\theta=4$)} & \makecell{10 \\[-5pt] (all $\theta$)} & \makecell{10 \\[-5pt] (all $\theta$)} & \makecell{36.15 \\[-5pt] ($\theta=4$)} & \makecell{ 9.735 \\[-5pt] ($\theta=0.25$) } & \makecell{ 4.204 \\[-5pt] ($\theta=0.25$) } \\
         \hline
    \end{tabular}
    \caption{Optimization results obtained by different parameters $\theta=0.25,0.5,1,2,4,$ and 8}
    \label{tab:param_theta}
\end{table}

\begin{table}[h!]
    \centering
    \begin{tabular}{c|c|c|c|c|c|c|c}
         \hline
         \multicolumn{2}{c|}{Ex.} & \makecell{Time \\[-5pt] (sec.)} & \makecell{No. of \\[-5pt] outer loops} & \makecell{No. of \\[-5pt]serious steps} & Cost & \makecell{${\hat {\bar p}}_f$ \\[-5pt] ($\cdot 10^{-4}$) } & \makecell{ ${\hat p}_f$ \\[-5pt] ($\cdot 10^{-4}$) } \\
         \hline
         \multirow{ 2}{*}{ \makecell{ \\[-5pt] 1 } } & Highest & \makecell{0.40 \\[-5pt] ($\omega=5$)} & \makecell{10 \\[-5pt] ($\omega=1.2$)} & \makecell{7 \\[-5pt] ($\omega=2$)} & \makecell{2,743 \\[-5pt] ($\omega=1.2$)} & \makecell{ 9.985 \\[-5pt] (all $\omega$) } & \makecell{ 3.05 \\[-5pt] ($\omega=5$) } \\ \cline{2-8}
         & Lowest & \makecell{0.14 \\[-5pt] ($\omega=1.2$)} & \makecell{8 \\[-5pt] ($\omega=5$)} & \makecell{6 \\[-5pt] ($\omega=1.2$)} & \makecell{2,718 \\[-5pt] ($\omega=3$)} & \makecell{ 9.985 \\[-5pt] (all $\omega$) } & \makecell{ 2.678 \\[-5pt] ($\omega=2$) } \\
         \hline
         \multirow{ 2}{*}{ \makecell{ \\[-5pt] 2 } } & Highest & \makecell{42 \\[-5pt] ($\omega=1.2$)} & \makecell{6 \\[-5pt] ($\omega=1.5$)} & \makecell{5 \\[-5pt] ($\omega=1.5$)} & \makecell{29.29 \\[-5pt] ($\omega=3$)} & \makecell{ 9.860 \\[-5pt] ($\omega=1.2$) } & \makecell{ 3.779 \\[-5pt] ($\omega=1.2$) } \\ \cline{2-8}
         & Lowest & \makecell{19 \\[-5pt] ($\omega=2$)} & \makecell{2 \\[-5pt] ($\omega=1.2$)} & \makecell{2 \\[-5pt] ($\omega=1.2$)} & \makecell{28.61 \\[-5pt] ($\omega=1.2$)} & \makecell{ 8.684 \\[-5pt] ($\omega=1.5$) } & \makecell{ 3.353 \\[-5pt] ($\omega=3$) } \\
         \hline
         \multirow{ 2}{*}{ \makecell{ \\[-5pt] 3 } } & Highest & \makecell{5.13 \\[-5pt] ($\omega=5$)} & \makecell{10 \\[-5pt] (all $\omega$)} & \makecell{10 \\[-5pt] (all $\omega$)} & \makecell{38.39 \\[-5pt] ($\omega=5$)} & \makecell{ 9.860 \\[-5pt] ($\omega=2$) } & \makecell{ 4.330 \\[-5pt] ($\omega=2$) } \\ \cline{2-8}
         & Lowest & \makecell{1.4 \\[-5pt] ($\omega=1.2$)} & \makecell{10 \\[-5pt] (all $\omega$)} & \makecell{10 \\[-5pt] (all $\omega$)} & \makecell{36.06 \\[-5pt] ($\omega=1.5$)} & \makecell{ 9.585 \\[-5pt] ($\omega=1.5$) } & \makecell{ 4.054 \\[-5pt] ($\omega=1.5$) } \\
         \hline
    \end{tabular}
    \caption{Optimization results obtained by different parameters $\omega=1.2,1.5,2,3,$ and 5}
    \label{tab:param_omega}
\end{table}

\subsection{Target buffered failure probability}

To test its stability to the magnitude of target buffered failure probability (i.e., $\bar{p}_f^t$), the algorithm is tested with target values $1\cdot 10^{-2}$ and $1\cdot 10^{-4}$. The optimization results are summarized in Table \ref{tab:change_bpf}. In the table, comparisons are made in parentheses with the results in Table \ref{tab:ex_results} with target probability $1\cdot 10^{-3}$. As expected, in all examples, the cost of obtained solutions increases with a higher $\bar{p}_f^t$. It is noted that computing time and the number of iterations remain stable, while the estimated buffered failure probabilities remain close but lower than the target values. This demonstrates the robustness of the algorithm with respect to the level of $\bar{p}_f^t$. 

\begin{table}[h!]
    \centering
    \begin{tabular}{c|c|c|c|c|c|c|c}
         \Xhline{3\arrayrulewidth}
         \multicolumn{8}{l}{$\bar{p}_f^t = 0.1 \cdot 10^{-2}$} \\
         \Xhline{3\arrayrulewidth}
         Ex. & \makecell{Time \\[-5pt] (sec.)} & \makecell{No. of \\[-5pt] outer loops} & \makecell{No. of \\[-5pt]serious steps} & \makecell{Computed \\[-5pt] solution} & Cost & ${\hat {\bar p}}_f$ & ${\hat p}_f$ \\
         \hline
         1 & \makecell{0.131 \\[-5pt] ($-$0.00527)} & \makecell{8 \\[-5pt] (+1)} & \makecell{8 \\[-5pt] (+1)} & (1092, 150.0) & \makecell{2,334 \\[-5pt] ($-$409.3)} & \makecell{ 9.975 \\[-5pt] $\cdot 10^{-3}$ } & \makecell{ 3.030 \\[-5pt] $\cdot 10^{-3}$ } \\
         \hline
         2 & \makecell{40.4 \\[-5pt] ($+26.8$) } & \makecell{2 \\[-5pt] ($-1$)} & \makecell{2 \\[-5pt] ($-1$)} & \makecell{(1.501, 1.090, \\[-5pt] 1.338, 1.000)} & \makecell{27.73 \\[-5pt] ($-$0.6300)} & \makecell{ 9.015 \\[-5pt] $\cdot 10^{-3}$ } & \makecell{ 3.636 \\[-5pt] $\cdot 10^{-3}$ } \\
         \hline
         3 & \makecell{3.21 \\[-5pt] (+0.387) } & \makecell{6 \\[-5pt] ($-4$)} & \makecell{6 \\[-5pt] ($-4$)} & \makecell{(6.435, 6.435, \\[-5pt] 6.435, 6.435, \\[-5pt] 2.218, 6.435)} & \makecell{34.39 \\[-5pt] ($-$1.809)} & \makecell{ 9.949 \\[-5pt] $\cdot 10^{-3}$ } & \makecell{ 4.268 \\[-5pt] $\cdot 10^{-3}$ } \\
         \Xhline{3\arrayrulewidth}
         \multicolumn{8}{l}{$\bar{p}_f^t = 0.1 \cdot 10^{-4}$} \\
         \Xhline{3\arrayrulewidth}
         1 & \makecell{0.561 \\[-5pt] (+0.425)} & \makecell{7 \\[-5pt] (+0)} & \makecell{7 \\[-5pt] (+0)} & (1471, 150.0) & \makecell{3,091 \\[-5pt] (+348.2)} & \makecell{ 9.976 \\[-5pt] $\cdot 10^{-5}$ } & \makecell{ 3.100 \\[-5pt] $\cdot 10^{-5}$ } \\
         \hline
         2 & \makecell{32.8 \\[-5pt] (+19.2) } & \makecell{3 \\[-5pt] (+0)} & \makecell{3 \\[-5pt] (+0)} & \makecell{(1.668, 1.000, \\[-5pt] 1.765, 1.000)} & \makecell{29.72 \\[-5pt] (+1.090)} & \makecell{ 9.801 \\[-5pt] $\cdot 10^{-5}$ } & \makecell{ 3.300 \\[-5pt] $\cdot 10^{-5}$ } \\
         \hline
         3 & \makecell{5.80 \\[-5pt] (+2.98) } & \makecell{8 \\[-5pt] ($-2$)} & \makecell{8 \\[-5pt] ($-2$)} & \makecell{(7.615, 7.616, \\[-5pt] 7.617, 7.616, \\[-5pt] 1.000, 7.616)} & \makecell{39.08 \\[-5pt] (+2.881)} & \makecell{ 9.801 \\[-5pt] $\cdot 10^{-5}$ } & \makecell{ 4.350 \\[-5pt] $\cdot 10^{-5}$ }  \\
         \hline
    \end{tabular}
    \caption{Optimization results obtained by different target values of buffered failure probability (in parentheses, comparisons are made with results in Table \ref{tab:ex_results}.)}
    \label{tab:change_bpf}
\end{table}

\subsection{Initial points}\label{sec:init_point}

Since problems of interest are nonconvex, the quality of computed solutions is most dependent on initial points. To test this, for each example, 100 experiments are performed with different initial points sampled by Latin hypercube sampling. The best solutions of each example, i.e., the solution that satisfies reliability constraints and incur the lowest cost, are summarized in Table \ref{tab:init_sols}. By considering a solution the best if it yields a cost less than 3\% higher than the lowest obtained value, it is observed that 100\%, 46\%, and 7\% of the initial points lead to one of the best solutions in the first, second, and third examples, respectively. This suggests that the third example has in particular multiple local optima. It is noted that the results in Table \ref{tab:ex_results}, which are obtained with initial solutions as a midpoint of the lower and upper bounds, are all one of the best obtained solutions. Nonetheless, the result strongly indicates that one needs to try multiple initializations to ensure the quality of an obtained solution. The proposed algorithm is advantageous to this end owing to its computational efficiency as demonstrated by various tests illustrated above.

\begin{table}[h!]
    \centering
    \begin{tabular}{c|c|c|c|c|c}
         \hline
         Ex. & \makecell{Best optimal \\[-5pt] solution$^\mathrm{a}$} & Cost & \makecell{${\hat {\bar p}}_f$ \\[-5pt] $(\cdot 10^{-4})$} & \makecell{${\hat p}_f$ \\[-5pt] $(\cdot 10^{-4})$} & \makecell{Ratio of best \\[-5pt] solutions$^\mathrm{b}$} \\
         \hline
         1 & (1257, 150.0) & 2,709 & 9.960 & 3.203 & 100 \% \\
         \hline
         2 & \makecell{(1.581, 1.000, \\[-5pt] 1.443, 1.000)} & 28.52 & 9.835 & 3.679 & 46 \% \\
         \hline
         3 & \makecell{(7.215, 6.896, \\[-5pt] 7.391, 7.391, \\[-5pt] 1.000, 6.906)} & 36.80 & 6.857 & 3.053 & 7 \% \\
         \hline
         \multicolumn{6}{l}{\small $^\mathrm{a}$ The feasible solution leading to the lowest cost} \\[-6pt]
         \multicolumn{6}{l}{\small $^\mathrm{b}$ Feasible solutions leading to a cost less than 3\% higher than the lowest one}
    \end{tabular}
    \caption{Optimization results with varying initial solutions}
    \label{tab:init_sols}
\end{table}

\section{Conclusions}

This study proposes an efficient algorithm for reliability-based optimization (RBO), particularly for general system events that are represented as a link-set of cut-sets. To handle such general systems, we evaluate system reliability by realizations of random variables (i.e., samples or data points) instead of analytical calculation that requires problem-specific formulas. This can be done by employing the buffered optimization and reliability method (BORM) that replaces the conventional failure probability by the buffered failure probability. We call it the {\it S-BORM} algorithm.

The S-BORM algorithm efficiently solves RBO problems of general systems by leveraging four ideas. First, a reliability constraint is penalized and moved to the objective function. Second, limit-state functions are linearized adaptively at a current solution. Third, the modified objective function is reformulated as a difference-of-convex function so that its optimization can be solved by a difference-of-convex bundle method. Fourth, an active-set strategy is employed, which enables us to consider only a small subset of samples that are within or close enough to failure domains. We do not assume convexity either for cost function or for limit-state functions. This makes the S-BORM algorithm applicable for a wide class of systems. Although such minimal assumptions make it difficult (if not impossible) to theoretically guarantee global optimality and convergence, we provide justifications of the proposed approach and empirical demonstrations by presenting newly designed numerical examples. Examples include complex and large-scale systems such as a truss structure system with 108 failure modes and an electrical system combined with reliability growth models. We show that the S-BORM algorithm provides a handy means to solve these complex problems. Moreover, extensive parametric investigations show that the algorithm remains insensitive to algorithm parameter values and magnitude of target failure probability.

Utilizing realizations of random variables greatly facilitates general applications as it eliminates the need for problem-specific derivations and enables non-parametric analysis. By leveraging such advantages, we develop a Matlab-based tool, which is available at \url{https://github.com/jieunbyun/sborm}. To run the algorithm, users only need to provide general information, i.e., cost functions and limit-state functions (including their gradients) and realizations of random variables. Meanwhile, this underlines a distinct potential of the BORM for further development of general software tools of reliability-based optimization. This is advantageous considering barriers of implementing RBO algorithms, which often require a high level of  knowledge and engineering. Promising topics for being combined with such data-driven optimization include surrogate models, real-time data, and sequential decision-making.

\section*{Acknowledgement}

The research of the first author is in part supported by the Humboldt Research Fellowship for Postdoctoral Researchers from Alexander von Humboldt Foundation. The second author acknowledges financial support from the Gaspard-Monge Program for Optimization and Operations Research (PGMO) project
“SOLEM - Scalable Optimization for Learning and Energy Management”. The research of the third author is supported in part by the Air Force Office of Scientific Research, Mathematical Optimization (21RT0484).

\section*{Conflicts of interest}

The authors have no relevant financial or non-financial interests to disclose.

\section*{Replication of results}

The examples can be replicated by (Matlab-based) codes and data uploaded at \url{https://github.com/jieunbyun/sborm}.

\bibliographystyle{elsarticle-num-names}
\bibliography{refs.bib}

\begin{thebibliography}{30}
\expandafter\ifx\csname natexlab\endcsname\relax\def\natexlab#1{#1}\fi
\providecommand{\url}[1]{\texttt{#1}}
\providecommand{\href}[2]{#2}
\providecommand{\path}[1]{#1}
\providecommand{\DOIprefix}{doi:}
\providecommand{\ArXivprefix}{arXiv:}
\providecommand{\URLprefix}{URL: }
\providecommand{\Pubmedprefix}{pmid:}
\providecommand{\doi}[1]{\href{http://dx.doi.org/#1}{\path{#1}}}
\providecommand{\Pubmed}[1]{\href{pmid:#1}{\path{#1}}}
\providecommand{\bibinfo}[2]{#2}
\ifx\xfnm\relax \def\xfnm[#1]{\unskip,\space#1}\fi
\bibitem[{Enevoldsen and Sørensen(1994)}]{rbo94}
\bibinfo{author}{I.~Enevoldsen}, \bibinfo{author}{J.~D. Sørensen},
\newblock \bibinfo{title}{Reliability-based optimization in structural
  engineering},
\newblock \bibinfo{journal}{Structural Safety} \bibinfo{volume}{15}
  (\bibinfo{year}{1994}) \bibinfo{pages}{169--196}.
\bibitem[{Byun and Song(2020)}]{mbnDm20}
\bibinfo{author}{J.-E. Byun}, \bibinfo{author}{J.~Song},
\newblock \bibinfo{title}{Efficient probabilistic multi-objective optimization
  of complex systems using matrix-based bayesian network},
\newblock \bibinfo{journal}{Reliability Engineering and System Safety}
  \bibinfo{volume}{200} (\bibinfo{year}{2020}) \bibinfo{pages}{106899}.
\bibitem[{Byun and Song(2021)}]{bns_jt21}
\bibinfo{author}{J.-E. Byun}, \bibinfo{author}{J.~Song},
\newblock \bibinfo{title}{A general framework of bayesian network for system
  reliability analysis using junction tree},
\newblock \bibinfo{journal}{Reliability Engineering and System Safety}
  \bibinfo{volume}{216} (\bibinfo{year}{2021}) \bibinfo{pages}{107952}.
\bibitem[{Lim and Song(2012)}]{sra12}
\bibinfo{author}{H.-W. Lim}, \bibinfo{author}{J.~Song},
\newblock \bibinfo{title}{Efficient risk assessment of lifeline networks under
  spatially correlated ground motions using selective recursive decomposition
  algorithm},
\newblock \bibinfo{journal}{Earthquake Engineering and Structural Dynamics}
  \bibinfo{volume}{41} (\bibinfo{year}{2012}) \bibinfo{pages}{1861--1882}.
\bibitem[{Nguyen et~al.(2011)Nguyen, Song, and Paulino}]{sl-rbo11}
\bibinfo{author}{T.~H. Nguyen}, \bibinfo{author}{J.~Song},
  \bibinfo{author}{G.~H. Paulino},
\newblock \bibinfo{title}{Single-loop system reliability-based topology
  optimization considering statistical dependence between limit-states},
\newblock \bibinfo{journal}{Structural and Multidisciplinary Optimization}
  \bibinfo{volume}{44} (\bibinfo{year}{2011}) \bibinfo{pages}{593--611}.
\bibitem[{Royset et~al.(2001)Royset, Der~Kiureghian, and Polak}]{rbo-decoup01}
\bibinfo{author}{J.~O. Royset}, \bibinfo{author}{A.~Der~Kiureghian},
  \bibinfo{author}{E.~Polak},
\newblock \bibinfo{title}{Reliability-based optimal structural design by the
  decoupling approach},
\newblock \bibinfo{journal}{Reliability Engineering and System Safety}
  \bibinfo{volume}{73} (\bibinfo{year}{2001}) \bibinfo{pages}{213--221}.
\bibitem[{Bismut et~al.(2022)Bismut, Pandey, and Straub}]{inm_pipe22}
\bibinfo{author}{E.~Bismut}, \bibinfo{author}{M.~D. Pandey},
  \bibinfo{author}{D.~Straub},
\newblock \bibinfo{title}{Reliability-based inspection and maintenance planning
  of a nuclear feeder piping system},
\newblock \bibinfo{journal}{Reliability Engineering and System Safety}
  \bibinfo{volume}{224} (\bibinfo{year}{2022}) \bibinfo{pages}{108521}.
\bibitem[{Kim and Song(2021)}]{rbo-kriging21}
\bibinfo{author}{J.~Kim}, \bibinfo{author}{J.~Song},
\newblock \bibinfo{title}{Reliability-based design optimization using quantile
  surrogates by adaptive gaussian process},
\newblock \bibinfo{journal}{ASCE Journal of Engineering Mechanics}
  \bibinfo{volume}{147} (\bibinfo{year}{2021}) \bibinfo{pages}{04021020}.
\bibitem[{Li et~al.(2019)Li, Gong, Gu, Jing, Fang, and Gao}]{rbo_surro_mcs19}
\bibinfo{author}{X.~Li}, \bibinfo{author}{C.~Gong}, \bibinfo{author}{L.~Gu},
  \bibinfo{author}{Z.~Jing}, \bibinfo{author}{H.~Fang},
  \bibinfo{author}{R.~Gao},
\newblock \bibinfo{title}{A reliability-based optimization method using
  sequential surrogate model and monte carlo simulation},
\newblock \bibinfo{journal}{Structural and Multidisciplinary Optimization}
  \bibinfo{volume}{59} (\bibinfo{year}{2019}) \bibinfo{pages}{439–--460}.
\bibitem[{Perdikaris et~al.(2015)Perdikaris, Venturi, Royset, and
  Karniadakis}]{Perdikaris}
\bibinfo{author}{P.~Perdikaris}, \bibinfo{author}{D.~Venturi},
  \bibinfo{author}{J.~O. Royset}, \bibinfo{author}{G.~E. Karniadakis},
\newblock \bibinfo{title}{Multi-fidelity modelling via recursive co-kriging and
  gaussian–markov random fields},
\newblock \bibinfo{journal}{Proc. R. Soc. A.} \bibinfo{volume}{471}
  (\bibinfo{year}{2015}) \bibinfo{pages}{20150018}.
\bibitem[{Canelas et~al.(2019)Canelas, Carrasco, and López}]{rbo_polytop19}
\bibinfo{author}{A.~Canelas}, \bibinfo{author}{M.~Carrasco},
  \bibinfo{author}{J.~López},
\newblock \bibinfo{title}{A new method for reliability analysis and
  reliability-based design optimization},
\newblock \bibinfo{journal}{Structural and Multidisciplinary Optimization}
  \bibinfo{volume}{59} (\bibinfo{year}{2019}) \bibinfo{pages}{1655--1671}.
\bibitem[{Ding et~al.(2021)Ding, Hu, and Li}]{rbo-ga21}
\bibinfo{author}{Y.~Ding}, \bibinfo{author}{Y.~Hu}, \bibinfo{author}{D.~Li},
\newblock \bibinfo{title}{Redundancy optimization for multi-performance
  multi-state series-parallel systems considering reliability requirements},
\newblock \bibinfo{journal}{Reliability Engineering and System Safety}
  \bibinfo{volume}{215} (\bibinfo{year}{2021}) \bibinfo{pages}{107873}.
\bibitem[{Li et~al.(2022)Li, Chi, and Yu}]{rbo-pso22}
\bibinfo{author}{S.~Li}, \bibinfo{author}{X.~Chi}, \bibinfo{author}{B.~Yu},
\newblock \bibinfo{title}{An improved particle swarm optimization algorithm for
  the reliability–redundancy allocation problem with global reliability},
\newblock \bibinfo{journal}{Reliability Engineering and System Safety}
  \bibinfo{volume}{225} (\bibinfo{year}{2022}) \bibinfo{pages}{108604}.
\bibitem[{Royset and Wets(2021)}]{optim_primer}
\bibinfo{author}{J.~O. Royset}, \bibinfo{author}{R.~J.-B. Wets},
  \bibinfo{title}{An Optimization Primer}, \bibinfo{edition}{1} ed.,
  \bibinfo{publisher}{Springer Cham}, \bibinfo{year}{2021}.
\bibitem[{Royset and Polak(2007)}]{rbo_samp07}
\bibinfo{author}{J.~O. Royset}, \bibinfo{author}{E.~Polak},
\newblock \bibinfo{title}{Extensions of stochastic optimization results to
  problems with system failure probability functions},
\newblock \bibinfo{journal}{Journal of Optimization Theory and Applications}
  \bibinfo{volume}{133} (\bibinfo{year}{2007}) \bibinfo{pages}{1--18}.
\bibitem[{Rockafellar and Royset(2010)}]{bpf10}
\bibinfo{author}{R.~T. Rockafellar}, \bibinfo{author}{J.~O. Royset},
\newblock \bibinfo{title}{On buffered failure probability in design and
  optimization of structures},
\newblock \bibinfo{journal}{Reliability Engineering and System Safety}
  \bibinfo{volume}{95} (\bibinfo{year}{2010}) \bibinfo{pages}{499--510}.
\bibitem[{Byun and Royset(2022)}]{borm22}
\bibinfo{author}{J.-E. Byun}, \bibinfo{author}{J.~O. Royset},
\newblock \bibinfo{title}{Data-driven optimization of reliability using
  buffered failure probability},
\newblock \bibinfo{journal}{Structural Safety} \bibinfo{volume}{98}
  (\bibinfo{year}{2022}) \bibinfo{pages}{102232}.
\bibitem[{Royset and Byun(2021)}]{bpf_grad21}
\bibinfo{author}{J.~O. Royset}, \bibinfo{author}{J.-E. Byun},
\newblock \bibinfo{title}{Gradients and subgradients of buffered failure
  probability},
\newblock \bibinfo{journal}{Operations Research Letters} \bibinfo{volume}{49}
  (\bibinfo{year}{2021}) \bibinfo{pages}{868--873}.
\bibitem[{Chaudhuri et~al.(2022)Chaudhuri, Kramer, Norton, Royset, and
  Willcox}]{rbdo_cvar22}
\bibinfo{author}{A.~Chaudhuri}, \bibinfo{author}{B.~Kramer},
  \bibinfo{author}{M.~Norton}, \bibinfo{author}{J.~O. Royset},
  \bibinfo{author}{K.~Willcox},
\newblock \bibinfo{title}{Certifiable risk-based engineering design
  optimization},
\newblock \bibinfo{journal}{AIAA Journal} \bibinfo{volume}{60}
  (\bibinfo{year}{2022}) \bibinfo{pages}{551--565}.
\bibitem[{Tuy(2016)}]{dc_textbook16}
\bibinfo{author}{H.~Tuy}, \bibinfo{title}{Convex Analysis and Global
  Optimization}, \bibinfo{edition}{2} ed., \bibinfo{publisher}{Springer},
  \bibinfo{year}{2016}.
\bibitem[{{de Oliveira}(2019)}]{pb17}
\bibinfo{author}{W.~{de Oliveira}},
\newblock \bibinfo{title}{Proximal bundle methods for nonsmooth {DC}
  programming},
\newblock \bibinfo{journal}{Journal of Global Optimization}
  \bibinfo{volume}{75} (\bibinfo{year}{2019}) \bibinfo{pages}{523--563}.
\bibitem[{{van Ackooij} and {de Oliveira}(2014)}]{vanAckooij_Oliveira_2014}
\bibinfo{author}{W.~{van Ackooij}}, \bibinfo{author}{W.~{de Oliveira}},
\newblock \bibinfo{title}{Level bundle methods for constrained convex
  optimization with various oracles},
\newblock \bibinfo{journal}{Computational Optimization and Applications}
  \bibinfo{volume}{57} (\bibinfo{year}{2014}) \bibinfo{pages}{555--597}.
\bibitem[{Pietrantuono et~al.(2010)Pietrantuono, Russo, and Trivedi}]{rgm10}
\bibinfo{author}{R.~Pietrantuono}, \bibinfo{author}{S.~Russo},
  \bibinfo{author}{K.~S. Trivedi},
\newblock \bibinfo{title}{Software reliability and testing time allocation: An
  architecture-based approach},
\newblock \bibinfo{journal}{IEEE Transactions on Software Engineering}
  \bibinfo{volume}{36} (\bibinfo{year}{2010}) \bibinfo{pages}{323--337}.
\bibitem[{Bertolino et~al.(2019)Bertolino, Miranda, Pietrantuono, and
  Russo}]{rgm19}
\bibinfo{author}{A.~Bertolino}, \bibinfo{author}{B.~Miranda},
  \bibinfo{author}{R.~Pietrantuono}, \bibinfo{author}{S.~Russo},
\newblock \bibinfo{title}{Adaptive test case allocation, selection and
  generation using coverage spectrum and operational profile},
\newblock \bibinfo{journal}{IEEE Transactions on Software Engineering}
  \bibinfo{volume}{47} (\bibinfo{year}{2019}) \bibinfo{pages}{881--898}.
\bibitem[{Song and Der~Kureghian(2003)}]{LpBound03}
\bibinfo{author}{J.~Song}, \bibinfo{author}{A.~Der~Kureghian},
\newblock \bibinfo{title}{Bounds on system reliability by linear programming},
\newblock \bibinfo{journal}{ASCE Journal of Engineering Mechanics}
  \bibinfo{volume}{129} (\bibinfo{year}{2003}) \bibinfo{pages}{627--636}.
\bibitem[{Der~Kureghian et~al.(2007)Der~Kureghian, Ditlevsen, and
  Song}]{power08}
\bibinfo{author}{A.~Der~Kureghian}, \bibinfo{author}{O.~D. Ditlevsen},
  \bibinfo{author}{J.~Song},
\newblock \bibinfo{title}{Availability, reliability and downtime of systems
  with repairable components},
\newblock \bibinfo{journal}{Reliability Engineering and System Safety}
  \bibinfo{volume}{92} (\bibinfo{year}{2007}) \bibinfo{pages}{231--242}.
\bibitem[{Byun and Song(2017)}]{MsrRgm17}
\bibinfo{author}{J.-E. Byun}, \bibinfo{author}{J.~Song},
\newblock \bibinfo{title}{Reliability growth analysis of k-out-of-n systems
  using matrix-based system reliability method},
\newblock \bibinfo{journal}{Reliability Engineering and System Safety}
  \bibinfo{volume}{165} (\bibinfo{year}{2017}) \bibinfo{pages}{410--421}.
\bibitem[{{de Oliveira}(2020)}]{abc-dc}
\bibinfo{author}{W.~{de Oliveira}},
\newblock \bibinfo{title}{The {ABC of DC} programming},
\newblock \bibinfo{journal}{Set-Valued and Variational Analysis}
  \bibinfo{volume}{28} (\bibinfo{year}{2020}) \bibinfo{pages}{679--706}.
\bibitem[{Han and Mangasarian(1979)}]{Han_Mangasarian_1979}
\bibinfo{author}{S.~Han}, \bibinfo{author}{O.~L. Mangasarian},
\newblock \bibinfo{title}{Exact penalty functions in nonlinar programming},
\newblock \bibinfo{journal}{Mathematical Programming} \bibinfo{volume}{17}
  (\bibinfo{year}{1979}) \bibinfo{pages}{251--269}.
\bibitem[{Le~Thi et~al.(2012)Le~Thi, Pham~Dinh, and Ngai}]{LeThi2012}
\bibinfo{author}{H.~A. Le~Thi}, \bibinfo{author}{T.~Pham~Dinh},
  \bibinfo{author}{H.~V. Ngai},
\newblock \bibinfo{title}{Exact penalty and error bounds in {DC} programming},
\newblock \bibinfo{journal}{Journal of Global Optimization}
  \bibinfo{volume}{52} (\bibinfo{year}{2012}) \bibinfo{pages}{509--535}.

\end{thebibliography}

\appendix
\section{Derivation of difference-of-convex decomposition}\label{sec:dc_deriv}

We derive a difference-of-convex decomposition of the objective function (\ref{eq:obj_pen}), so that the optimization problem can be formulated as a difference-of-convex program. Recall that the function is defined with a current point $\hat{\bm{x}}^\nu$ and a penalization constant $\theta^\nu$ at an iteration step $\nu$. The last term of the function can be decomposed as
\begin{linenomath*}
\begin{equation}
\begin{aligned}
    \min_{q\in\mathds{Q}_k} [ \beta_q^\nu(\bm{v}) + \langle \alpha_q^\nu(\bm{v}), \bm{x} \rangle] -\gamma &= \sum_{q\in\mathds{Q}_k} [ \beta_q^\nu(\bm{v}) + \langle \alpha_q^\nu(\bm{v}), \bm{x} \rangle] -\gamma - \max_{q\in\mathds{Q}_k} \sum_{r\in\mathds{Q}_k\backslash\{q\}}[ \beta_r^\nu(\bm{v}) + \langle \alpha_r^\nu(\bm{v}), \bm{x} \rangle] \label{eq:decomp1} \\
    &:= -p_k^\nu(\bm{x},\gamma,\bm{v}).
\end{aligned}
\end{equation}
\end{linenomath*}
Since the first and second terms in (\ref{eq:decomp1}) are linear and convex functions, respectively, we have a convex function
\begin{linenomath*}
\begin{equation}\label{eq:pk_def}
    p_k^\nu(\bm{x},\gamma,\bm{v}) = \max_{q\in\mathds{Q}_k}\left\{ \sum_{r\in\mathds{Q}_k\backslash\{q\}}[ \beta_r^\nu(\bm{v}) + \langle \alpha_r^\nu(\bm{v}), \bm{x} \rangle] \right\}- \sum_{q\in\mathds{Q}_k}[ \beta_q^\nu(\bm{v}) + \langle \alpha_q^\nu(\bm{v}), \bm{x} \rangle] + \gamma.
\end{equation}
\end{linenomath*}
Let $q^*$ be an index yielding the maximum in (\ref{eq:pk_def}). Then, a subgradient of $p_k^\nu$ at $(\bm{x},\gamma)$ can be computed as
\begin{linenomath*}
\postdisplaypenalty=0
\begin{equation}
    g_{p_k}^\nu(\bm{v}) := \left( \begin{matrix} \sum_{r\in\mathds{Q}_k\backslash\{q^*\}} \alpha_r^\nu(\bm{v}) - \sum_{q\in\mathds{Q}_k} \alpha_q^\nu(\bm{v}) \\ 1  \end{matrix} \right) = \left( \begin{matrix} -\alpha_{q^*}^\nu(\bm{v}) \\ 1 \end{matrix} \right).
\end{equation}
\end{linenomath*}

We now move on to the next decomposition, i.e.,
\begin{linenomath*}
\postdisplaypenalty=0
\begin{align*}
    \max_{k=1,\ldots,K} \min_{q\in\mathds{Q}_k} \beta_q^\nu(\bm{v}) + \langle \alpha_q^\nu(\bm{v}),\bm{x} \rangle &= \max_{k=1,\ldots,K} -p_k^\nu(\bm{x},\gamma,\bm{v}) \\
    &= \max_{k=1,\ldots,K} \sum_{l \neq k}^K p_l^\nu(\bm{x},\gamma,\bm{v}) - \sum_{k=1}^K p_k^\nu(\bm{x},\gamma,\bm{v}) \\
    &= \psi^\nu(\bm{x},\gamma,\bm{v}) - \varphi^\nu(\bm{x},\gamma,\bm{v}),
\end{align*}
\end{linenomath*}
where both $\psi^\nu(\bm{x},\gamma,\bm{v})$ and $\varphi^\nu(\bm{x},\gamma,\bm{v})$ are convex. Let $k^*$ be an index yielding the maximum above.
Then,
\[
    g_{\varphi}^\nu(\bm{v}) := \left( \sum_{k=1}^K g_{p_k}^\nu (\bm{v}) \right) \in \partial \varphi^\nu(\bm{x},\gamma,\bm{v}), \quad \mbox{and}\quad
    g_\phi^\nu(\bm{v}) := [g_\varphi^\nu(\bm{v})-g_{p_{k^*}}^\nu(\bm{v})] \in \partial \psi^\nu(\bm{x},\gamma,\bm{v}).
\]
Since $\max\{0, \psi^\nu(\bm{x},\gamma,\bm{v}) - \varphi^\nu(\bm{x},\gamma,\bm{v})\}=\max\{\psi^\nu(\bm{x},\gamma,\bm{v}), \varphi^\nu(\bm{x},\gamma,\bm{v})\}-\varphi^\nu(\bm{x},\gamma,\bm{v})$, we can write
\begin{linenomath*}
\postdisplaypenalty=0
\begin{align*}
    &\gamma + \frac{1}{\bar{p}_f^t} \sum_{n\in\hat{\mathds{N}}^\nu} p_n \max\Big\{0, \max_{k\in{1,\ldots,K}} \min_{q\in\mathds{Q}_k} \beta_q^\nu(\bm{v}_n) + \langle \alpha_q^\nu(\bm{v}_n), \bm{x} \rangle -\gamma \Big\} \\
    &= \gamma + \frac{1}{\bar{p}_f^t} \sum_{n\in\hat{\mathds{N}}^\nu} p_n \Big\{ \max \big\{ \psi^\nu(\bm{x},\gamma,\bm{v}_n), \varphi^\nu(\bm{x},\gamma,\bm{v}_n) \big\} - \varphi^\nu(\bm{x},\gamma,\bm{v}_n) \Big\} \\
    &= \gamma + \frac{1}{\bar{p}_f^t} \sum_{n\in\hat{\mathds{N}}^\nu} p_n \max \big\{ \psi^\nu(\bm{x},\gamma,\bm{v}_n), \varphi^\nu(\bm{x},\gamma,\bm{v}_n) \big\} - \frac{1}{\bar{p}_f^t} \sum_{n=1}^N p_n \varphi^\nu(\bm{x},\gamma,\bm{v}_n) \\
    &= \Gamma^\nu(\bm{x},\gamma) - \Lambda^\nu(\bm{x},\gamma),
\end{align*}
\end{linenomath*}
which decomposes into two convex functions
\begin{linenomath*}
\postdisplaypenalty=0
\begin{align*}
    \Gamma^\nu(\bm{x},\gamma) &= \gamma + \frac{1}{\bar{p}_f^t} \sum_{n\in\hat{\mathds{N}}^\nu} p_n \max \big\{ \psi^\nu(\bm{x},\gamma,\bm{v}_n), \varphi^\nu(\bm{x},\gamma,\bm{v}_n) \big\}, \\
    \Lambda^\nu(\bm{x},\gamma) &= \frac{1}{\bar{p}_f^t} \sum_{n\in\hat{\mathds{N}}^\nu} p_n \varphi^\nu(\bm{x},\gamma,\bm{v}_n).
\end{align*}
\end{linenomath*}
Their subgradients at $(\bm{x},\gamma)$ can be computed as
\begin{linenomath*}
\postdisplaypenalty=0
\begin{subequations}\label{eq:subgrads}
\begin{align}
    g_\Gamma^\nu &= \left( \begin{matrix} \mathbf{0}_{D\times1} \\ 1 \end{matrix} \right) + \frac{1}{\bar{p}_f^t}\Bigg[ \sum_{n\in\mathds{N}_x^\psi} p_n g_\psi^\nu(\bm{v}_n) + \sum_{n\in\hat{\mathds{N}}^\nu\backslash\mathds{N}_x^\psi} p_n g_\varphi^\nu(\bm{v}_n) \Bigg], \\
    g_\Lambda^\nu &= \frac{1}{\bar{p}_f^t} \Bigg[ \sum_{n=1}^N p_n g_\varphi^\nu(\bm{v}_n) \Bigg],
\end{align}
\end{subequations}
\end{linenomath*}
where $\mathds{N}_x^{\psi,\nu}=\{n\in\hat{\mathds{N}}^\nu: \psi^\nu(\bm{x},\gamma,\bm{v}_n)\geq\varphi^\nu(\bm{x},\gamma,\bm{v}_n)\}$.

Finally, recall that $c(\bm{x})$ has $L$-Lipschitz continuous gradients over $\mathds{X}$. Thus, $c(\bm{x})+ \frac{L}{2}\|\bm{x}\|^2$ is convex and yields  the following difference-of-convex decomposition for $c$: $c(\bm{x})= c(\bm{x})+ \frac{L}{2}\|\bm{x}\|^2 - \frac{L}{2}\|\bm{x}\|^2$. Furthermore, by writing $\max\{ 0, \Gamma^\nu(\bm{x},\gamma)-\Lambda^\nu(\bm{x},\gamma) \}=\max \{ \Gamma^\nu(\bm{x},\gamma), \Lambda^\nu(\bm{x},\gamma) \}-\Lambda^\nu(\bm{x},\gamma)$, the approximated objective function (\ref{eq:obj_pen}) becomes
\begin{linenomath*}
\begin{subequations}\label{dc-decomposotion}
\begin{equation}\label{eq:obj_pen_final}
\begin{split}
    F^\nu(\bm{x},\gamma;\theta^\nu,\hat{\bm{x}}^\nu) &= c(\bm{x}) + \theta^\nu \max \big\{ \Gamma^\nu(\bm{x},\gamma), \Lambda^\nu(\bm{x},\gamma) \big\}-\theta^\nu\Lambda^\nu(\bm{x},\gamma) \\
    &= c(\bm{x}) + \frac{L}{2} \norm{\bm{x}}^2+\theta^\nu \max \big\{ \Gamma^\nu(\bm{x},\gamma), \Lambda^\nu(\bm{x},\gamma) \big\} - \bigg[ \frac{L}{2} \norm{\bm{x}}^2 + \theta^\nu\Lambda^\nu(\bm{x},\gamma) \bigg] \\
    &= f_1^\nu(\bm{x},\gamma) - f_2^\nu(\bm{x},\gamma),
\end{split}
\end{equation}
where
\postdisplaypenalty=0
\begin{align}
    f_1^\nu(\bm{x},\gamma) &= c(\bm{x}) + \frac{L}{2} \norm{\bm{x}}^2+\theta^\nu \max \big\{ \Gamma^\nu(\bm{x},\gamma), \Lambda^\nu(\bm{x},\gamma) \big\}, \\
    f_2^\nu(\bm{x},\gamma) &=\frac{L}{2} \norm{\bm{x}}^2 + \theta^\nu\Lambda^\nu(\bm{x},\gamma),
\end{align}
\end{subequations}
\end{linenomath*}
and $f_1^\nu(\bm{x},\gamma)$ and $f_2^\nu(\bm{x},\gamma)$ are convex functions with a straightforward rule for computing a pair of subgradients.
Hence, with its objective as a difference-of-convex function, problem~\eqref{rbo:bpf_penal_linear} is a difference-of-convex program, which can be solved by Algorithm 1 in \cite{pb17}. 

The derived difference-of-convex decomposition uses the Lipschitz constant $L$ of the function $c$, a value that may be unknown to the decision-maker. However, as $c(\bm{x})+ \frac{\bar L}{2}\|\bm{x}\|^2$ is still convex for all $\bar L \geq L$ \cite[Prop.1]{abc-dc}, we only need a rough overestimation of this constant to obtain a difference-of-convex decomposition $c(\bm{x})+ \frac{\bar L}{2}\|\bm{x}\|^2 - \frac{\bar L}{2}\|\bm{x}\|^2$ for $c$, and thus (in view of \eqref{dc-decomposotion}) for the objective function of \eqref{rbo:bpf_penal_linear}. We mention that when $c$ is convex (the case in our numerical experiments), we can simply set $\bar{L}=0$ in \eqref{dc-decomposotion} and the resulting difference-of-convex decomposition is still valid.

\section{Algorithmic Justifications}\label{appendix:dc}

The function $F^\nu(\bm{x},\gamma;\theta^\nu,\hat{\bm{x}}^\nu)$, appearing in \ref{sborm_al:pb_eval} of the algorithm, can be written as the difference of two convex functions as required by \cite[Alg. 1]{pb17}; the formula is derived in \ref{sec:dc_deriv}. For limit-state functions that are linear in $\bm{x}$, S-BORM can only converge to critical points as defined in \cite{pb17}. Further motivation for the algorithm follows below. 

Generally, $F^\nu(\bm{x},\gamma;\theta^\nu,\hat{\bm{x}}^\nu)$ is an approximation of the actual function $F(\bm{x},\gamma;\theta^\nu)$ but it becomes exact at the current point $\hat{\bm{x}}^\nu$, i.e.,
\begin{linenomath*}
\begin{equation*}
    F^\nu(\hat{\bm{x}}^\nu,\gamma;\theta^\nu,\hat{\bm{x}}^\nu)= F(\hat{\bm{x}}^\nu,\gamma ;\theta^\nu),\quad \forall\; \gamma.
\end{equation*}
\end{linenomath*}
Therefore, since $(\hat{\bm{x}}^\nu,\hat \gamma^{\nu})$ is feasible to the subproblem of \ref{sborm_al:pb_eval}, we have that
\begin{linenomath*}
\begin{equation*}
F^\nu(\bm{x}^{\nu+1},\gamma^{\nu+1};\theta^\nu,\hat{\bm{x}}^\nu)+ \frac{\lambda^\nu}{2} \|\bm{x}^{\nu+1} - \hat{\bm{x}}^\nu\|_2^2+\frac{\lambda^\nu}{2}(\gamma^{\nu+1} -\hat \gamma^\nu)^2\leq  F^{\nu}(\hat{\bm{x}}^{\nu},\hat \gamma^{\nu};\theta^\nu,\hat{\bm{x}}^\nu) = F(\hat{\bm{x}}^\nu,\hat \gamma^\nu;\theta^\nu).
\end{equation*}
\end{linenomath*}
This means that the predicted decrease $\zeta^\nu$ is non-negative and
S-BORM is a descent method in this sense.

If \ref{sborm_al:pb_eval} produces $(\bm{x}^{\nu+1},\gamma^{\nu+1}) = (\hat{\bm{x}}^{\nu},\hat{\gamma}^{\nu})$, then one can show that this point is  critical for the penalized and linearized problem~\eqref{rbo:bpf_penal_linear} (with $\theta=\theta^\nu$ and $\hat{\bm{x}}=\hat{\bm{x}}^\nu$). A similar conclusion appears to hold more generally too. Suppose that $\tt{tol}=0$ and consider two cases. 

\paragraph{1. The algorithm produces only finitely many serious steps followed by an infinite sequence of null steps} Then, $(\hat{\bm{x}}^\nu,\hat{\gamma}^\nu)$ equals some fixed point $(\hat{\bm{x}},\hat{\gamma})$ for all $\nu$ after the last serious step. Furthermore, the linearization of limit-state functions are fixed for these iterations. The S-BORM algorithm therefore reduces to a penalized approach for solving \begin{linenomath*}
\begin{subequations}\label{rbo:bpf_penal_linear_wlo}
\postdisplaypenalty=0
\begin{alignat}{3}
& \!\min_{\bm{x}\in\mathds{X},\gamma\in \mathds{R}} &\qquad & c(\bm{x}) \\
& \text{subject to} &      & \gamma + \frac{1}{\bar{p}_f^t} \sum_{n=1}^N p_n \max\Big\{0, \max_{k\in{1,\ldots,K}} \min_{q\in\mathds{Q}_k} g_q(\hat{\bm{x}},\bm{v}_n) + \langle \nabla g_q(\hat{\bm{x}},\bm{v}_n), \bm{x}-\hat{\bm{x}} \rangle -\gamma \Big\} \leq 0.
\end{alignat}
\end{subequations}
\end{linenomath*}
As the penalized parameter becomes $\theta^{\max}$ after finitely many steps, the iterates produced by \ref{sborm_al:pb_eval}  converges to $(\hat{\bm{x}},\hat{\gamma})$ as $\lambda^\nu$ increases indefinitely and the value  $F^\nu(\bm{x},\gamma;\theta^{\max},\hat{\bm{x}})$
given in \eqref{eq:obj_pen} does not change  (because both the  penalization parameter and the current solution candidate are fixed for all $\nu$ large enough). This argument suggests that $(\hat{\bm{x}},\hat{\gamma})$ is a critical point for the linearized problem. The theory of exact penalty functions asserts that if a condition on the constraint of \eqref{rbo:bpf_penal_linear_wlo} exists, and $\theta^{\max}$ is greater than the largest optimal dual variable
associated with this constraint, then solutions of the penalized problem
\eqref{rbo:bpf_penal_linear} are solutions of \eqref{rbo:bpf_penal_linear_wlo}; see, e.g., \cite{Han_Mangasarian_1979,LeThi2012} and \cite[Proposition 6.13]{optim_primer}. 

\paragraph{2. The algorithm produces infinitely many serious steps}
Again, recall that the penalty parameter becomes $\theta^{\max}$ after finitely many iterations.
In this case, as already argued, the descent test ensures that the sequence of function values is non-increasing:
\begin{linenomath*}
\[
F(\hat{\bm{x}}^{\nu+1},\hat \gamma^{\nu+1};\theta^{\max}) \leq F(\hat{\bm{x}}^\nu,\hat \gamma^\nu;\theta^{\max}) -\kappa \zeta^\nu, \quad \mbox{with $\zeta^\nu$ given in \ref{sborm_al:eval_sol}.}
\]
\end{linenomath*}
With the mild assumption that $F(\cdot,\cdot;\theta^{\max})$ has bounded level sets, a simple recursive argument (the telescope sum) on the above inequality shows that $\zeta^\nu \to 0$ and one can argue that any cluster point of the sequence of serious steps $\{(\hat{\bm{x}}^\nu,\hat \gamma^\nu)\}$ is critical for the penalized problem~\eqref{rbo:bpf_penal_linear}. Once again, the penalization arguments concerning \eqref{rbo:bpf_penal_linear} and \eqref{rbo:bpf_penal_linear_wlo} apply.

\paragraph{Summary}
The above arguments tell us that S-BORM always terminates after finitely many steps provided that $\tt{tol}>0$. When $\tt{tol}=0$ and the limit-state functions are linear in $\bm{x}$, then the algorithm is guaranteed to converge to a critical point. For more general limit-state functions, the discussion on the asymptotic behavior of the algorithm suggests that the approach computes a critical point for the linearized problem \eqref{rbo:bpf_penal_linear_wlo}, with $\hat{\bm{x}}$ being the $x$-part of the last serious step, or an arbitrary cluster point (if any) of $\{\hat{\bm{x}}^\nu\}_{\nu=1}^\infty$. A full mathematical argument is beyond the scope of the present paper.



\end{document}